\documentclass[12pt]{amsart}
\usepackage{amsfonts}
\usepackage{mathrsfs}
\usepackage{amscd,amsmath,amssymb,amsfonts}
\usepackage{color}
\usepackage[all]{xy}
\theoremstyle{plain}
\newtheorem{thm}{Theorem}
\newtheorem{lem}[thm]{Lemma}
\newtheorem{cor}[thm]{Corollary}
\newtheorem{prop}[thm]{Proposition}

\theoremstyle{definition}
\newtheorem{defn}[thm]{Definition}

\newtheorem{rmk}[thm]{Remark}

\numberwithin{thm}{section} \numberwithin{equation}{section}

\newcommand{\ga}[2]{\begin{gather}\label{#1}#2 \end{gather}}

\newcommand{\sB}{{\mathcal B}}

\newcommand{\sD}{{\mathcal D}}
\newcommand{\sE}{{\mathcal E}}
\newcommand{\sF}{{\mathcal F}}

\newcommand{\sH}{{\mathcal H}}
\newcommand{\sI}{{\mathcal I}}

\newcommand{\sK}{{\mathcal K}}
\newcommand{\sL}{{\mathcal L}}
\newcommand{\sM}{{\mathcal M}}

\newcommand{\sO}{{\mathcal O}}
\newcommand{\sP}{{\mathcal P}}
\newcommand{\sQ}{{\mathcal Q}}
\newcommand{\sR}{{\mathcal R}}

\newcommand{\sU}{{\mathcal U}}
\newcommand{\sV}{{\mathcal V}}
\newcommand{\sW}{{\mathcal W}}

\def\wt#1{\widetilde {#1}}

\begin{document}

\title{Frobenius splitting of moduli spaces of parabolic bundles}
\author{Xiaotao Sun}
\address{Center of Applied Mathematics, School of Mathematics, Tianjin University, No.92 Weijin Road, Tianjin 300072, P. R. China}
\email{xiaotaosun@tju.edu.cn}
\author{Mingshuo Zhou}
\address{Center of Applied Mathematics, School of Mathematics, Tianjin University, No.92 Weijin Road, Tianjin 300072, P. R. China}
\email{zhoumingshuo@amss.ac.cn}
\date{May 2, 2023}
\thanks{Both authors are supported by the NSFC (No.11831013, No.11921001 and No.12171352).}
\begin{abstract} Let $C$ be a nonsingular projective curve over an algebraically closed field of characteristic $p>0$ and $I\subset C$ be a finite set. If
$\mathcal{U}_{C,\,\omega}$ denotes the moduli space of semistable parabolic bundles of rank $r$ and degree $d$ on $C$ with parabolic structures determined by
$\omega=(k,\{\vec n(x),\vec a(x)\}_{x\in I})$, we prove that $\mathcal{U}_{C,\,\omega}$ is \textit{$F$-split} for generic $C$ and generic choice of $I$ when $p>3r$.
\end{abstract}
\keywords{Frobenius split,  Moduli spaces, Parabolic bundles}
\subjclass{Algebraic Geometry, 14H60, 14D20}
\maketitle
\begin{quote}

\end{quote}
\section{Introduction}

Let $M$ be a variety over a perfect field of characteristic $p>0$ and $F:M\to M$ be the absolute Frobenius morphism. Then $M$ is called
\textit{$F$-split }(or \textit{split}) if $\sO_M\hookrightarrow F_*\sO_M$ splits. If there exists a $\sigma\in {\rm H}^0(\omega^{-1}_M)$ such that
$\sigma^{p-1}$ provides a splitting of ${\sO_M\hookrightarrow F_*\sO_M}$, we call that $M$ is\textit{ split by a $(p-1)$-power}.
Although some important varieties are proved to be $F$-split, for example, flag varieties and their Schubert subvarieties (cf. \cite{MeRa}, \cite{RR}), the product of
two flag varieties for the same group $G$ (cf. \cite{MeRa88}) and cotangent bundles of flag varieties (cf. \cite{KLT}), we do not know many other examples.
For example, the question if  moduli space of semistable bundles on a curve $C$ is $F$-split remains open.

For a bundle $E$ on $C$, parabolic structures at $x\in C$ are given by flags $E_x=Q_{l_x+1}(E)_x\twoheadrightarrow
Q_{l_x}(E)_x\twoheadrightarrow\cdots\twoheadrightarrow
Q_1(E)_x\twoheadrightarrow Q_0(E)_x=0$ of fixed type $\vec n(x)=(n_1(x),\ldots,n_{l_x+1}(x))$ and a sequence of integers
$$\vec a(x)=(a_1(x),\ldots,a_{l_x+1}(x)),\quad 0=\frac{a_1(x)}{k}<\frac{a_2(x)}{k}<\cdots<\frac{a_{l_x+1}(x)}{k}<1,$$
where $n_i(x):={\rm dim}({\rm ker}\{Q_i(E)_x\twoheadrightarrow Q_{i-1}(E)_x\})$. For any finite set $I$ of points on $C$ and $r\ge 2$, moduli space $\sU_{C,\,\omega}$ of  semistable parabolic bundles of rank $r$ and degree $d$ on $C$, with parabolic structures determined by $\omega:=(k,\{\vec n(x),\vec a(x)\}_{x\in I})$,
is a normal projective variety.

When $r=2$ and $p>5$, Mehta-Ramadas proved in \cite{M.R} that $\sU_{C,\,\omega}$ is $F$-split for generic curve $C$ and generic parabolic points on $C$.
In this article, we generalize the result to higher rank of parabolic bundles.

\begin{thm}\label{thm1.1} Let $h\ge 2g$ be a positive integer and $r\ge 2$. Then, for a generic curve $C$ of genus $g$ over an algebraically closed field of characteristic $p>3r$ and
a set $I=\{y_1,\ldots,y_h\}\cup\{z_1,\ldots,z_{2+h-2g}\}$ of generic points on $C$ with $\vec n(y_i)=(r-1,1)$, $\vec n(z_i)=(1,\ldots,1)$, the moduli space $\sU_{C,\,\omega_c}$
of semistable parabolic bundles of rank $r$ and degree $d$ on $C$, with parabolic structures determined by canonical weight $\omega_c$, is split by a $(p-1)$-power.
\end{thm}

\begin{cor}\label{cor1.2} For any $r\ge 2$, when $p>3r$, $\sU_{C,\,\omega}$ is $F$-split for generic $C$ and generic parabolic points on $C$.
\end{cor}

To see how Theorem \ref{thm1.1} implies  Corollary \ref{cor1.2}, for any $m\ge 0$ and $I'=\{z_1,\ldots,z_m\}$, let $\omega=(k,\{\vec n'(x),\vec a(x)\}_{x\in I'}$ be
arbitrary type and weights at $x\in I'$. Then the projections ${\rm Flag}_{(1,\ldots,1)}(\sE_{z_i})\to {\rm Flag}_{\vec n'(z_i)}(\sE_{z_i})$
induce a morphism $\sU_{C,\,\omega_c}\supset U\xrightarrow{f} \sU_{C,\,\omega}$ such that $f_*\sO_U=\sO_{\sU_{C,\,\omega}}$ when $h$ is large enough ($g$ is fixed)
(see Lemma \ref{lem2.8}). Thus Corollary \ref{cor1.2} is true for $\sU_{C,\,\omega}$ if Theorem \ref{thm1.1} holds for $\sU_{C,\,\omega_c}$.

To explain the proof of Theorem \ref{thm1.1}, let $(\wt X, \wt I)$ be a smooth $|\wt I|$-pointed curve of genus $\wt{g}=g-1$, $\wt\omega_c=(2r,\{\vec n(x),\vec a_c(x)\}_{x\in \wt I})$, and
assume that there is a $\wt\sigma\in {\rm H}^0(\omega^{-1}_{\sU_{\wt X,\,\wt \omega_c}})$ such that $\wt\sigma^{p-1}$ splits $\sU_{\wt X,\,\wt \omega_c}$. Then, take $x_1,\,x_2 \in \wt I$ and identify them to get a $(|\wt I|-2)$-pointed curve $(X, I)$ (where $I=\wt I\setminus\{x_1,\,x_2\}\subset X$) of genus $g$ with only one node $x_0$, we will produce, start with $\wt\sigma\in {\rm H}^0(\omega^{-1}_{\sU_{\wt X,\,\wt \omega_c}})$, a $\sigma\in {\rm H}^0(\Theta_{\sU_{X,\,\omega_c}})$ such that $\sigma^{p-1}$ splits an open set of $\sU_{X,\,\omega_c}$ where $\omega_c=(2r, \{\vec n(x),\,\vec a_c(x)\}_{x\in I})$. Then, take a flat family of pointed curves $C_T\to T$ with sections $T_i\subset C_T$ ($1\le i\le |I|$) such that $C_0=X$, $C_t$ ($t\neq 0$) are smooth curves of genus $g$ and
$\{\,T_i\cap C_0\,|\,1\le i\le |I|\,\}=I\subset X$, we will construct a $T$-flat scheme $\sU_{C_T,\,\omega_c}\xrightarrow{\nu} T$ and a line bundle $\Theta_{\sU_{C_T,\,\omega_c}}$ on $\sU_{C_T,\,\omega_c}$ such that
\ga{1.1} {\nu^{-1}(t)=\sU_{C_t,\,\omega_c},\quad \Theta_{\sU_{C_T,\,\omega_c}}|_{\sU_{C_t,\,\omega_c}}=\Theta_{\sU_{C_t,\,\omega_c}}\quad (\forall\,\,t\in T),}
where we assume that there is an $x \in I$ with $\vec n(x)=(1,\ldots,1)$ so that our GIT quotient commutes with base changes (see Lemma \ref{lem2.13} and Remark \ref{rmk2.14}).
Then, by ${\rm H}^1(\Theta_{\sU_{X,\,\omega_c}})=0$, the section $\sigma\in{\rm H}^0(\Theta_{\sU_{X,\,\omega_c}})$ spreads to fibers $\sU_{C_t,\,\omega_c}$ nearby $\sU_{X,\,\omega_c}$, which splits $\sU_{C_t,\,\omega_c}$ by a $(p-1)$-power by Proposition \ref{prop2.2}.

We remark that key step in above arguments is the construction of $\sigma\in {\rm H}^0(\Theta_{\sU_{X,\,\omega_c}})$ from $\wt\sigma\in {\rm H}^0(\omega^{-1}_{\sU_{\wt X,\,\wt \omega_c}})$, which we need essentially $\vec n(x_1)=\vec n(x_2)=(1,\ldots,1)$. Thus, in order to proceed with our arguments, we have to assume
that $\wt I$ contains at least $3$ points of type $(1,\ldots,1)$ (we need an $x\in I=\wt I\setminus\{x_1,x_2\}$ with $\vec n(x)=(1,\ldots,1)$ so that \eqref{1.1} holds). Hence, to complete the proof of our main theorem, we have to show that $\mathcal{U}_{\mathbb{P}^1,\,{\omega}_c}$ is split by a $(p-1)$-power for $\omega_c=(2r, \{\vec n(x),\,\vec a_c(x)\}_{x\in I})$ such that $I$ contains \textit{many points $x\in I$ with $\vec n(x)=(1,\ldots,1)$}. Unfortunately, even for $I=\{z_1,z_2,z_3\}$ with $\vec n(z_i)=(1,\ldots,1)$, it is not immediately clear that $\mathcal{U}_{\mathbb{P}^1,\,{\omega}_c}$ is split by a $(p-1)$-power when $r>2$.
One of the technical results in this article is to prove that $\mathcal{U}_{\mathbb{P}^1,\,{\omega}_c}$ is split by a $(p-1)$-power for $I=\{y_1, z_1,z_2,z_3\}$ with
$\vec n(y_1)=(r-1,1)$ and $\vec n(z_i)=(1,\ldots,1)$ (see Theorem \ref{thm5.1}). Then, by a similar degeneration argument, we have

\begin{thm}\label{thm1.3} For any $h>0$, let $I\subset \mathbb{P}^1$ be a set of $2+2h$ points with
$\vec n(y_1)=\cdots=\vec n(y_{h})=(r-1,1)$ $(y_i\in I)$ and  $\vec n(x)=(1,\ldots, 1)$ $(\forall x\in I-\{y_1,\ldots, y_h\})$.
Then, when $p>3r$, the moduli space $\mathcal{U}_{\mathbb{P}^1,\,{\omega}_c}$ is split by a $(p-1)$-power for generic choice of points in $I$.
\end{thm}

We will prove it by induction on $h$. Assume theorem holds for $h_1=1$ and $h_2=h-1$, let $X_1=X_2=\mathbb{P}^1$ and $I_i\subset X_i$ (with $|I_i|=h_i$) such that
Theorem \ref{thm1.3} holds for $(X_i,I_i)$. Take point $x_i\in I_i$ with $\vec n(x_i)=(1,\ldots,1)$ and glue $X_1$, $X_2$ by identifying $x_1$ and $x_2$:
$$\wt X:=X_1\sqcup X_2\xrightarrow{\pi} X:=X_1\cup X_2,\quad \pi^{-1}(x_0)=\{x_1,\,x_2\},$$
where $\{x_0\}=X_1\cap X_2$ is the only node of $X$, we will construct
$$\sigma\in {\rm H}^0(\Theta_{\sU_{X,\,\omega_c}}), \quad \omega_c=(2r,\{\vec n(x),\vec a_c(x)\}_{x\in (I_1\setminus\{x_1\})\cup(I_2\setminus\{x_2\})}),$$
from $\sigma_i\in {\rm H}^0(\omega^{-1}_{\sU_{X_i,\,\omega^i_{c}}})$, where $\omega^i_c=(2r,\{\vec n(x),\vec a_c(x)\}_{x\in I_i})$ and $\sigma_i^{p-1}$ splits
$\sU_{X_i,\,\omega^i_{c}}$ ($i=1,\,2$). Take a flat family $C_T\to T$ of $|I|$-pointed $\mathbb{P}^1$ with sections $T_i\subset C_T$ such that
$C_t=\mathbb{P}^1$ ($t\neq 0$), $C_0=X$ and
$$\{\,T_i\cap X\,|\,1\le i\le |I|\,\}=(I_1\setminus\{x_1\})\cup(I_2\setminus\{x_2\}):=I,$$
by a similar degeneration argument, $\sigma^{p-1}$ induces $F$-splitting of fibers $\sU_{C_t,\omega^t_c}$ nearby $\sU_{X,\,\omega_c}$, where $\omega^t_c=(2r,\{\vec n(x),\vec a_c(x)\}_{x\in I^t})$ with $$I^t=\{\,T_i\cap C_t\,|\,1\le i\le |I|\,\},\,\vec n(T_i\cap C_t)=\vec n(T_i\cap C_0).$$
Note $|I^t|=|I_1|-1+|I_2|-1=2+2h$ and $I^t$ contains $h+2$ points of type $(1,\ldots,1)$, Theorem \ref{thm1.3} is proved.

We describe briefly the contents of this article. In Section 2, some preliminary results are collected. In particular, we make an important observation
(see Remark \ref{rmk2.14}) which seems to be new. In Section 3, under the assumptions that Theorem \ref{thm3.8}, Theorem \ref{thm3.17} and the case $h=1$ of Theorem \ref{thm1.3} hold, we prove Theorem \ref{thm1.1} by degeneration arguments. In Section 4, Theorem \ref{thm3.8} and Theorem \ref{thm3.17} are proved. In Section 5, we prove the case $h=1$ of Theorem \ref{thm1.3} (see Theorem \ref{thm5.1}), which is trivial when $r=2$ (see Remark \ref{rmk5.5}).

\section{Preliminaries and variation of moduli spaces}

Let $M$ be a variety over an algebraically closed field of characteristic $p>0$ and $F: M\rightarrow
M$ be the absolute Frobenius morphism, which is the identity on the
underlying topology space $M$ and $p$-th power map on the
structure sheaf $\sO_M$.
\begin{defn}\label{defn2.1} $M$ is called \emph{Frobenius split} ($F$-split or split simply) if $F^{\sharp}:\sO_M\rightarrow F_{\ast}\sO_M$ splits, i.e., there exists an $\sO_M$-linear map $$\varphi:F_{\ast}\sO_M\rightarrow \sO_M\quad \text{such that}\quad \varphi \circ F^{\sharp}={\rm id}_{\sO_M}.$$
\end{defn}
When $M$ is smooth, we have $\mathrm{Hom}(F_{\ast}\sO_M,\sO_M)={\rm H}^0(\omega_M^{1-p})$.
Given a $\sigma\in {\rm H}^0(\omega_M^{-1})$, let $T_{\sigma}:F_{\ast}\sO_M\rightarrow \sO_M$ be homomorphism corresponding to $\sigma^{p-1}$. If $T_{\sigma}$ splits $\sO_M\rightarrow F_{\ast}\sO_M$, we call that $M$ is \emph{split by a $(p-1)$-power}.
A normal variety $M$ is called to be \emph{split by a $(p-1)$-power} if open set $M^{sm.}$ of smooth points is so.
We summarize some facts about $F$-splitting (see \cite{M.S} and \cite[Proposition 2.1]{M.R}).

\begin{prop}\label{prop2.2}
(1) Let $f: M\rightarrow N$ be a morphism with
$f_{\ast}\sO_M=\sO_N$. Then $N$ is $F$-split if $M$ is so.

(2) If $M$ is normal, $U\subset M$ is open and ${\rm Codim}(M\setminus U)\ge 2$, then $M$ is split (resp. {split by a $(p-1)$-power}) if and only if $U$ is so.

(3) If $M$ is a proper F-split variety, then ${\rm H}^i(M,L)=0$ $(i\geq 1)$ for any ample line bundle $L$ on $M$.

(4) Let $M_S\rightarrow S$ be a flat morphism with $S$ smooth and generic closed fiber is normal with only constant global functions. Given a section of the relative anti-canonical bundle, if this section splits a nonempty open set of a fiber, then it splits fibers in a neighborhood.

(5) Let $M_S\rightarrow S$ be a smooth morphism with $S$ smooth (resp. generic closed fiber has only constant global sections). Given a section of anti-canonical bundle, if this section splits $M_S$, its restriction splits a nonempty open subset of generic fiber (resp. splits generic fiber).
\end{prop}

Let $C$ be an irreducible projective curve over an algebraically closed field of characteristic $p>0$ with at most one node. Let $I$ be a finite set of smooth points of $C$, and $E$ be a torsion free sheaf of rank $r$ and degree $d$ on $C$. The rank $r(E)$ is defined to be dimension of $E_{\xi}$ at generic point $\xi\in C$, and $d=\chi(E)-r(1-g)$.
\begin{defn}\label{defn2.3} By a quasi-parabolic structure of $E$ at a
smooth point $x\in C$, we mean a choice of flag of quotients
$$E_x=Q_{l_x+1}(E)_x\twoheadrightarrow
Q_{l_x}(E)_x\twoheadrightarrow\cdots\cdots\twoheadrightarrow
Q_1(E)_x\twoheadrightarrow Q_0(E)_x=0$$ of fiber $E_x$, $n_i(x)={\rm
dim}(\mathrm{ker}\{Q_i(E)_x\twoheadrightarrow Q_{i-1}(E)_x\})$ ($1\le i\le l_x+1$)
are called type of the flags. If, in addition, a sequence of integers
$$0= a_1(x)<a_2(x)<\cdots
<a_{l_x+1}(x)< k$$ are given, we call that $E$ has a parabolic
structure of type $$\vec n(x)=(n_1(x),n_2(x),\cdots,n_{l_x+1}(x))$$ and
weight $\vec a(x)=(a_1(x),a_2(x),\cdots,a_{l_x+1}(x))$ at $x\in C$.
\end{defn}

\begin{defn}\label{defn2.4} For subsheaf $F\subset E$, let $Q_i(E)_x^F\subset
Q_i(E)_x$ be the image of $F$ and $n_i^F(x)={\rm
dim}(\mathrm{ker}\{Q_i(E)_x^F\twoheadrightarrow Q_{i-1}(E)_x^F\})$. Let
$${\rm par}\chi(E):=\chi(E)+\frac{1}{k}\sum_{x\in
I}\sum^{l_x+1}_{i=1}a_i(x)n_i(x),$$
$${\rm par}\chi(F):=\chi(F)+\frac{1}{k}\sum_{x\in
I}\sum^{l_x+1}_{i=1}a_i(x)n^F_i(x).$$
Then $E$ is called \emph{semistable} (resp. \emph{stable}) for $\omega=(k, \{\vec n(x),\vec a(x)\}_{x\in I})$ if for any
nontrivial $F\subset E$ such that $E/F$ is torsion free,
one has
$${\rm par}\chi(F)\leq
\frac{{\rm par}\chi(E)}{r(E)}\cdot r(F)\,\,(\text{resp. }<).$$
The weight $\omega$ is called generic, if the semi-stability of parabolic bundles is equivalent to stability respect to $\omega=(k, \{\vec n(x),\vec a(x)\}_{x\in I})$.
\end{defn}

\begin{thm}[see Theorem 2.13 of \cite{Su3}]\label{thm2.5} There exists a seminormal projective variety $\mathcal{U}_{C,\, \omega}:=\mathcal{U}_{C}(r,d,\{k, \vec n(x), \vec a(x)\}_{x\in I}),$ which is the coarse moduli space of $s$-equivalent classes of semistable parabolic sheaves $E$ of rank $r$ and $\chi(E)=d+r(1-g)$ on $C$ with parabolic structures at $\{x\}_{x\in I}$ determined by $\omega$. When $C$ is smooth,  $\mathcal{U}_{C,\, \omega}$ is a normal projective variety.
\end{thm}

Recall the construction of $\sU_{C,\,\omega}$. Fix a line bundle $\sO(1)=\sO_C(c\cdot y)$ on $C$ of ${\rm deg}(\sO(1))=c$, let
$\chi=d+r(1-g)$, $P$ denote the polynomial $P(m)=crm+\chi$ and $V=k^{\oplus P(N)}$. Assume that $N$ is large enough such that $F(N)$ is generated by global sections and ${\rm H}^1(F(N))=0$ for any semistable parabolic sheaf $F$ of degree $d$. Let $\bold Q$ be the Quot scheme of quotients $V\otimes\sO_{C}(-N)\to F\to 0$ (of rank
$r$ and degree $d$) on $C$ and $V\otimes\sO_{C\times\bold Q}(-N)\to \sF\to 0$ be the universal quotient on $C\times\bold Q$.
Let $\sF_x=\sF|_{\{x\}\times\bold Q}$ and ${\rm Flag}_{\vec n(x)}(\sF_x)\to\bold Q$ be the relative flag scheme of type $\vec n(x)$. Then the group ${\rm SL}(V)$ acts on
$$\sR=\underset{x\in I}{\times_{\bold Q}}{\rm Flag}_{\vec n(x)}(\sF_x)\to \bold Q$$
and the data $\omega=(k, \{\vec n(x),\,\,\vec a(x)\}_{x\in I})$ (the
weight $(k,\{\vec a(x)\}_{x\in I})$ more precisely) determines a polarization
$$\Theta_{\sR,\omega}=({\rm
det}R\pi_{\sR}\sE)^{-k}\otimes\bigotimes_{x\in I}
\lbrace\bigotimes^{l_x}_{i=1} {\rm det}(\sQ_{\{x\}\times
\sR,i})^{d_i(x)}\rbrace\otimes\bigotimes_q{\rm
det}(\sE_y)^{\ell}$$
on $\sR$ such that the open set $\sR^{ss}_{\omega}$ (resp. $\sR^s_{\omega}$) of
GIT semistable (resp. GIT stable) points are precisely the set of semistable (resp. stable) parabolic sheaves on $C$ (see \cite{Su3}), where $\sE$ is the pullback of $\sF$ (under
 $C\times\sR\to C\times \bold Q$), ${\rm
det}R\pi_{\sR}\sE$ is determinant line bundle of cohomology,
$$\sE_x=\sQ_{\{x\}\times \sR,l_x+1}\twoheadrightarrow\sQ_{\{x\}\times \sR,l_x}\twoheadrightarrow \sQ_{\{x\}\times \sR,l_x-1}
\twoheadrightarrow\cdots\twoheadrightarrow \sQ_{\{x\}\times
\sR,1}\twoheadrightarrow0$$ are universal flags on $\sR$ of type $\vec n(x)$, $d_i(x)=a_{i+1}(x)-a_i(x)$, $r_i(x)=n_1(x)+n_2(x)+\cdots +n_i(x)$ and
$$\ell:=\frac{k\chi-\sum_{x\in I}\sum^{l_x}_{i=1}d_i(x)r_i(x)}{r}.$$
Then $\sU_{C,\,\omega}$ (resp. $\mathcal{U}^s_{C,\, \omega}$) is the GIT quotient $\sR^{ss}_{\omega}\xrightarrow{\psi} \sR^{ss}_{\omega}//{\rm SL}(V):=\sU_{C,\,\omega}$ (resp. $\sR^{s}_{\omega}\xrightarrow{\psi}\sR^{s}_{\omega}//{\rm SL}(V):=\sU^s_{C,\,\omega}$). If $\ell$ is an integer, we have

\begin{thm} \label{thm2.6} Assume $|I|\neq 0$ and $\vec n(x)=(1,\ldots, 1)$ for some $x\in I$. Then there exists a unique ample line bundle $\Theta_{\mathcal{U}_{C,\, \omega}}$ on
$\mathcal{U}_{C,\, \omega}$ such that $\psi^{\ast}\Theta_{\mathcal{U}_{C,\,\omega}}=\Theta_{\sR^{ss}_{\omega}}$  where $\sR^{ss}_{\omega}\xrightarrow{\psi} \sU_{C,\,\omega}$ is the quotient map and
$$\Theta_{\sR^{ss}_{\omega}}=({\rm
det}R\pi_{\sR}\sE)^{-k}\otimes\bigotimes_{x\in I}
\lbrace\bigotimes^{l_x}_{i=1} {\rm det}(\sQ_{\{x\}\times
\sR,i})^{d_i(x)}\rbrace\otimes\bigotimes_y{\rm
det}(\sE_y)^{\ell}$$
for a fixed smooth point $y\in C$.
\end{thm}
\begin{proof} The only difference between the case of characteristic $p>0$ and the case of characteristic zero is that the descent lemma does not hold in positive characteristic. But, by \cite[Proposition 2.2]{M.R}, both Luna's {\'e}tale slice theorem
and descent lemma still hold in positive characteristic when stabilizers of points $(E, \{Q_{\bullet}(E)_x\}_{x\in I})\in \sR^{ss}_{\omega},$ which have closed orbits, are linearly reductive.

For point $(E, \{Q_{\bullet}(E)_x\}_{x\in I})\in \sR^{ss}_{\omega}$ with closed orbit, its stabilizer must be linearly reductive if there is a point $x\in I$ such that the flag $Q_{\bullet}(E)_x$ at $x\in I$ has type $\vec n(x)=(1,\ldots, 1).$ Indeed, if
$$(E, \{Q_{\bullet}(E)_x\}_{x\in I})=(E^{(1)}, \{Q^{(1)}_{\bullet}(E)_x\}_{x\in I})\oplus \cdots \oplus (E^{(s)}, \{Q^{(s)}_{\bullet}(E)_x\}_{x\in I})$$
is a direct sum of stable parabolic bundles, then for any $i\neq j$ the two stable parabolic bundles $(E^{(i)}, \{Q^{(i)}_{\bullet}(E)_x\}_{x\in I})$ and
$(E^{(j)}, \{Q^{(j)}_{\bullet}(E)_x\}_{x\in I})$ are not isomorphic to each other.
\end{proof}
Let $C$ be a smooth curve of genus $g$, and $\mathcal{R}_F\subset \mathcal{R}$ be the open subset consisting of locally free sheaves over $C$. Then, by the proof of \cite[Proposition 5.1]{S1}, we have
$$\mathrm{Codim}(\sR_F\setminus \sR^{ss}_{\omega})> \underset{1<r_1<r}{\mathrm{min}}\left\{r_1(r-r_1)(g-1)+\sum_{x\in I}\Sigma_{x,r_1}(\omega)\right\},$$ $\mathrm{Codim}(\sR^{ss}_{\omega}\setminus \sR^{s}_{\omega})\geq  \underset{1<r_1<r}{\mathrm{min}}\left\{r_1(r-r_1)(g-1)+\underset{x\in I}{\sum}\,\Sigma_{x,r_1}(\omega)\right\}$, where
\ga{2.1}{\Sigma_{x,r_1}(\omega):=\left\{\aligned &\sum^{l_x+1}_{j=1}\left(r_1-\sum^j_{i=1}m_i(x)\right)(n_j(x)-m_j(x))\\
&+\sum^{l_x+1}_{j=1}\left(r_1n_j(x)-rm_j(x)\right)\frac{a_j(x)}{k}\endaligned\right\}}
$0\leq m_i(x)\leq n_i(x)$ and $r_1=m_1(x)+\cdots +m_{l_x+1}(x)$. By \cite[Lemma 5.2]{S1}, $\Sigma_{x,r_1}(\omega)\ge\frac{1}{k}$ holds for
all $0<r_1<r$. Thus we have

\begin{prop}\label{prop2.7} Let $C$ be a smooth curve of genus $g$, and $\mathcal{R}_F\subset \mathcal{R}$ be the open subset consisting of locally free sheaves over $C$. Then, for any data $\omega=(k, \{\vec n(x), \vec a(x)\}_{x\in I})$, we have
\begin{itemize}
 \item[(1)]  $\,\,{\rm Codim}(\sR_{\omega}^{ss}\setminus \sR_{\omega}^s)\ge \underset{0<r_1<r}{\mathrm{min}}\{r_1(r-r_1)(g-1)\}+\frac{1}{k}|{\rm
 I}|$;
\item[(2)]  $\,\,{\rm Codim} (\sR_F\setminus\sR_{\omega}^{ss})>\underset{0<r_1<r}{\mathrm{min}}\{r_1(r-r_1)(g-1)\}+\frac{1}{k}|{\rm
I}|$.
\end{itemize}
\end{prop}

\begin{lem}\label{lem2.8} Let $C$ be a smooth curve of genus $g$, for subset $I'\subset C$ and $\omega'=(k',\{\vec n'(x),\vec a'(x)\}_{x\in I'})$ such that $I\subset I'$, $\vec n'(x)=(1,\ldots,1)$ when $x\in I$, and $\underset{0<r_1<r}{\mathrm{min}}\{r_1(r-r_1)(g-1)\}+\frac{|I'|}{k'}>1.$
Then, when $\sU_{C,\,\omega'}$ is $F$-split, $\sU_{C,\,\omega}$ is $F$-split.
\end{lem}

\begin{proof} Let $\mathcal{R}'=\underset{x\in I'}{\times_{\bold{Q}}}{\rm Flag}_{\vec n'(x)}(\mathcal{F}_x)\xrightarrow{\hat{f}} \mathcal{R}=\underset{x\in I}{\times_{\bold{Q}}}{\rm Flag}_{\vec n(x)}(\mathcal{F}_x)$ be the projection, which induces a morphism $f:U\rightarrow \mathcal{U}_{C,\, \omega}$ where
$$U:=\hat{U}//\mathrm{SL}(V)\subset \sU_{C,\,\omega'},\quad \hat{U}={\hat{f}}^{-1}(\mathcal{R}^{ss}_{\omega})\cap {\mathcal{R}'}^s_{\omega'}
\subset {\mathcal{R}'}^s_{\omega'} $$
are open subsets of $\sU_{C,\,\omega'}$ and ${\mathcal{R}'}^s_{\omega'}$. Then
it is easy to prove that $f_*\sO_U=\sO_{\sU_{C,\,\omega}}$.
In fact, for any open subset $V\subset \mathcal{U}_{C,\, \omega}$,
we have $${\rm H}^0(V,\sO_{\mathcal{U}_{C,\,\omega}})={\rm H}^0(\psi^{-1}(V), \sO_{\mathcal{R}^{ss}_{\omega}})^{inv.}={\rm H}^0({\hat{f}}^{-1}\psi^{-1}(V),\sO_{{\hat{f}}^{-1}(\sR^{ss}_{\omega})})^{inv.}$$
since $\hat{f}_*\sO_{{\hat{f}}^{-1}(\sR^{ss}_{\omega})}=\sO_{ \sR^{ss}_{\omega}}$. We have
$\mathrm{Codim}({\hat{f}}^{-1}(\mathcal{R}^{ss}_{\omega})\setminus \hat{U})\geq 2$ by Proposition \ref{prop2.7}. Thus $\sO_{\sU_{C,\,\omega}}(V)={\rm H}^0({\hat{f}}^{-1}\psi^{-1}(V),\sO_{\hat{U}})^{inv.}=\sO_U(f^{-1}(V))$. We are done by Proposition \ref{prop2.2}.
\end{proof}

\begin{rmk}\label{rmk2.9} In order to prove that $\sU_{C,\,\omega}$ is $F$-split, we can assume
\begin{itemize}
\item [(1)] There are points $x\in I$ with fixed special type $\vec n(x)$ we want, for example, the points $x\in I$ with type
$\vec n(x)=(1,...,1)$.
\item [(2)] $\omega=(2r,\{\vec n(x), \vec a(x)\}_{x\in I})$ satisfies $a_1(x)=0$ and
$$a_{i+1}(x)-a_i(x)=n_i(x)+n_{i+1}(x)\quad (1\leq i\leq l_x).$$
\item [(3)] The number $|I|$ can be large enough. For example, it satisfies
$$\underset{0<r_1<r}{\mathrm{min}}\{r_1(r-r_1)(g-1)+\sum_{x\in I}\Sigma_{x,r_1}(\omega_c)\}>1.$$

\end{itemize}
\end{rmk}
\begin{rmk}\label{rmk2.10} The weight $\omega$ in Remark \ref{rmk2.9} will be used throughout this paper, which we denote by
$\omega_c=(2r,\{\vec n(x), \vec a_c(x)\}_{x\in I})$. A precise formula of $\Sigma_{x,r_1}(\omega_c)$ for type $\vec n(x)=(r-1,1)$ is given by
$$\Sigma_{x,r_1}(\omega_c)= \begin{cases} (r-r_1)-(r-r_1)\frac{1}{2}=\frac{r-r_1}{2}, &\text{if $m_2(x)\neq 0$}\\
(r_1-r\cdot 0)\frac{1}{2}=\frac{r_1}{2},&\text{if $m_2(x)=0$  }.\end{cases}$$
When $\vec n(x)=(1,\ldots,1)$, let $r\ge j_{r_1}>j_{r_1-1}>\cdots>j_2>j_1\ge 1$ be the numbers such that $m_{j_k}(x)\neq 0$ for any $j_k\in S:=\{j_1,\,j_2,\,...,\,j_{r_1}\}$ and
$m_j(x)=0$ if $j\notin S$. Then we have (note $m_j(x)$ can only be $0$ or $1$)
$$\aligned &\Sigma_{x,r_1}(\omega_c)=\sum_{j\notin S}\left(r_1-\sum^j_{i=1}m_i(x)\right)+\sum_{j\notin S}\frac{r_1(j-1)}{r}+\sum_{j\in S}\frac{(r_1-r)(j-1)}{r}\\
&=\sum_{j\notin S}r_1+\sum^r_{j=1}\frac{r_1(j-1)}{r}-\sum_{j\notin S}\sum^j_{i=1}m_i(x)-\sum_{j\in S}(j-1)\\
&=\frac{1}{2}r_1(3r-2r_1-1)-\sum_{j\notin S}\sharp\{j_k\in S\,|\,j_k<j\,\}-\sum^{r_1}_{k=1}(j_k-1)\\
&=\frac{1}{2}r_1(3r-2r_1-1)-r_1(r-j_{r_1})-\sum^{r_1-1}_{k=1}k(j_{k+1}-j_k-1)-\sum^{r_1}_{k=1}(j_k-1)\\
&=\frac{r_1(r-r_1)}{2}+(r_1-1)j_{r_1}-\sum^{r_1-1}_{k=1}(kj_{k+1}-(k-1)j_k)=\frac{r_1(r-r_1)}{2}.
\endaligned$$
We call the generic weight $\omega=(k,\{\vec n(x),\vec{a}'(x)\}_{x\in I})$ is a minor modification of $\omega_c$ if $|\frac{a'_i(x)}{k}-\frac{\bar{a}_i(x)}{2r}|>0$ is {\em sufficiently small} ($1\leq i\leq l_x+1$) for all $x\in I$. In this case, we have $\sR^{ss}_{\omega}=\sR^s_{\omega}$ and
$$\aligned      \mathrm{Codim}(\sR_F\setminus \sR^{s}_{\omega})\geq \underset{1<r_1<r}{\mathrm{min}}\{r_1(r-r_1)(g-1)+\Sigma_{x\in I}\Sigma_{x,r_1}(\omega_c)\}.      \endaligned$$
\end{rmk}

Let $\mathrm{Det}:\sR_F\rightarrow J_C^d$
be the determinant map, which induces a morphism $\mathrm{Det}:\sU_{C,\,\omega}\rightarrow J_C^d$. Let $\sL$ be the universal line bundle on $C\times J_C^d$, and $\Theta_y\equiv(\mathrm{det}R\pi_J \sL)\otimes \sL_y^{d+1-g}.$ Then we have

\begin{lem}\label{lem2.11} Let $\omega^{-1}_{\sR_F}$ be the anti-canonical bundle of $\sR_F$. If $$\omega_c=(2r,\{\vec n(x), \vec a_c(x)\}_{x\in I})$$ satisfies the condition in Remark \ref{rmk2.9} (3). Then $$\omega^{-1}_{\sR_F}=\Theta_{\sR_F,\,\omega_c}\otimes (\mathrm{Det}^{\ast}\Theta_y)^{-2},\quad \omega^{-1}_{\sU_C,\,\omega_c}=\Theta_{\sU_C,\,\omega_c}\otimes (\mathrm{Det}^{\ast}\Theta_y)^{-2}.$$
\end{lem}
\begin{proof} This is in fact a reformulation of \cite[Proposition 2.2]{S1} where anti-canonical bundle $\omega^{-1}_{\sR_F}$
is given. The descent lemma holds since $\vec n(x)=(1,\ldots,1)$ for some $x\in I$ (see the proof of Theorem \ref{thm2.6}). Note that the quotient $\sR^s_{\omega_c}\rightarrow \mathcal{U}^s_{C,\,\omega_c}$ is a principal $\mathrm{PGL}(V)$-bundle, and $\mathrm{Codim}(\mathcal{U}_{C,\,\omega_c}\setminus \mathcal{U}_{C,\,\omega_c}^{s})\geq 2$ holds.
\end{proof}

\begin{defn}\label{defn2.12} The weight $\omega_c=(2r,\{\vec n(x), \vec a_c(x)\}_{x\in I})$ in Remark \ref{rmk2.9} (2) is called \textbf{canonical weight} determined by $\{\vec n(x)\}_{x\in I}$.
\end{defn}

Let $G$ be a reductive group acting on a variety $X$ such that the good quotient $\pi:X\rightarrow X//G$ (respect to a polarization) exists. Let $Y\subset X$ be a $G$-invariant closed subvariety, and $\widetilde{Y}$ be its (reduced) image in $X//G$ such that $Y=\pi^{-1}(\widetilde{Y})$. In general, the morphism
$\theta:Y//G\rightarrow \widetilde{Y}$ is only a (set-theoretically) bijective map (see \cite[Lemma A.1.2]{M}).
\begin{lem}\label{lem2.13} Assume that closed orbit in $X$ has linearly reductive isotropy group. Then $\theta:Y//G\rightarrow \widetilde{Y}$
is an isomorphism.
\end{lem}
\begin{proof} Take $y\in \widetilde{Y}\subset X//G$, let $x\in X$ such that $\pi(x)=y$ with closed orbit $G\cdot x$. By Luna's $\acute{e}\text{tale}$ slice theorem, there exists a scheme $T\rightarrow X$ whose image containing $x$ such that $T//G_x\rightarrow X//G$ is $\acute{e}\text{tale}$ (see \cite[Theorem 4.5]{A}).
By base change $\widetilde{Y}\rightarrow X//G$, we obtain
an $\acute{e}\text{tale}$ morphism $\varphi:T//G_x\times_{X//G} \widetilde{Y}\rightarrow \widetilde{Y}$.
Then in order to prove $\theta$ is an isomorphism, note it is bijective, we only need to show that $\varphi$ factors through $\theta$.
In fact, the reduced scheme $T//G_x\times_{X//G} \widetilde{Y}$ is
$$(T\times_{X//G}\widetilde{Y})//G_x=(T\times_{X//G}\widetilde{Y})^{\mathrm{red}}//G_x$$
and there is a factorization
$\varphi: (T\times_{X//G}\widetilde{Y})^{\mathrm{red}}//G_x \rightarrow Y//G\xrightarrow{\theta}\widetilde{Y} $
since $Y=\pi^{-1}(\widetilde{Y})$. We are done.\end{proof}

\begin{rmk}\label{rmk2.14} For $\omega=(k,\{\vec n(x), \vec a(x)\}_{x\in I})$, if there is a point $x\in I$ with type $\vec n(x)=(1,...,1)$, an important observation in
this paper is that any point $(E, \{Q_{\bullet}(E)_x\}_{x\in I})\in \sR^{ss}_{\omega}$ with closed orbit must have linearly reductive stabilizer (see the proof of Theorem \ref{thm2.6}).
Then, by \cite[Proposition 2.2]{M.R} and Lemma \ref{lem2.13}, both Luna's {\'e}tale slice theorem
and descent lemma hold in our case and the GIT quotient
$$\sU_{C,\omega}=\sR^{ss}_{\omega}//{\rm SL}(V)$$
commutes with base changes.
\end{rmk}

We recall some facts of representations of $\mathrm{GL}(r)$ in characteristic $p>0$.
Let $H\subset \mathrm{GL}(r)$ be subgroup of diagonal matrices and $P$ be the group of all characters of $H$. Then we have
$P=\mathbb{Z}\epsilon_1\oplus \mathbb{Z}\epsilon_2\oplus \cdots \oplus \mathbb{Z}\epsilon_r,$
where $\epsilon_i$ is a character defined by $\epsilon_i(\mathrm{diag}(\lambda_1,\ldots, \lambda_r))=\lambda_i$ ($1\leq i\leq r$).
Set $P^{\ast}=\mathrm{Hom}(P,\mathbb{Z})=\mathbb{Z}\epsilon^{\ast}_1\oplus \mathbb{Z}\epsilon^{\ast}_2\oplus \cdots \oplus \mathbb{Z}\epsilon^{\ast}_r$,  $h_i=\epsilon^{\ast}_i-\epsilon^{\ast}_{i+1}$ and $\alpha_i=\epsilon_i-\epsilon_{i+1}$ for $1\leq i\leq r-1$. Let $\alpha_0=\epsilon_1-\epsilon_r$, $h_0=\epsilon^{\ast}_1-\epsilon^{\ast}_r$ and
$$\aligned P^{+}=&\{\lambda\in P| \lambda(h_i)\geq 0 \,\, \text{for any}\,\, 1\leq i\leq r\}\\
C=&\{\lambda\in P^+| \lambda(h_0)\leq p-r+1\}.
\endaligned
$$
Given a $\lambda\in P^+$, in addition to irreducible representation $L(\lambda)$, there are other three representations associated with $\lambda$, namely the standard module $H^0(\lambda)$, the Weyl module $W(\lambda)$, and the tilting module $P(\lambda)$. If any two of $L(\lambda)$, $H^{0}(\lambda)$, $W(\lambda)$ and $P(\lambda)$ are isomorphic, then all four are isomorphic.
By Strong Linkage Principle,
for any $\lambda\in C$, Weyl module $W(\lambda)$ is irreducible and its dual is again a Weyl module.

\begin{lem}[see Proposition~10 of \cite{Ma}]\label{lem2.15} Let $\lambda\in C$ be a dominant weight and $\mathrm{char}({k})=p>r$. Then $P(\lambda)\cong W(\lambda)$ and $P(\lambda)$ is irreducible.
\end{lem}

\begin{lem}\label{lem2.16} When $p>r+m$, the representation $${\rm H}^0(Gr, \sO_{Gr}(m)):={\rm H}^0(\mathrm{Grass}_r(V_1\oplus V_2), (\mathrm{det}Q)^m)$$ of $\mathrm{GL}(V_1)\times \mathrm{GL}(V_2)$ is complete reducible.
\end{lem}

\begin{proof}
When $p>r+m$, ${\rm H}^0(Gr, \sO_{Gr}(m))$ is the irreducible representation $W_{\lambda}(V_1\oplus V_2)$ of $\mathrm{GL}(V_1\oplus V_2)$ with highest weight $\lambda=(m,\ldots, m)$ by Lemma \ref{lem2.15}. Moreover for any $\mu=(\mu_1,\ldots, \mu_r)$ with $0\leq \mu_r\leq \cdots \leq \mu_1\leq m$, the representations $W_{\mu}(V_1)$ and $W_{\nu}(V_2)$ of ${\rm GL}(r)$ are all irreducible where $\nu=(m-\mu_r,\ldots, m-\mu_1)$. Thus $W_{\mu}(V_1)\otimes W_{\nu}(V_2)$ are irreducible representations of $\mathrm{GL}(V_1)\times \mathrm{GL}(V_2)$ and
$$W_{\lambda}(V_1\oplus V_2)=\bigoplus_{\mu} W_{\mu}(V_1)\otimes W_{\nu}(V_2)$$
where $\mu=(\mu_1,\ldots, \mu_r)$ runs through all $0\leq \mu_r\leq \cdots \leq \mu_1\leq m$.
\end{proof}

\section{The degeneration arguments and normalization}

Let $\omega_c=(2r,\{\vec n(x), \vec a_c(x)\}_{x\in I})$ be the canonical weight and $\sU_{C,\,\omega_c}$ be the
moduli space of semistable parabolic bundles of rank $r$ and degree $d$ on $C$ with parabolic structures determined by $\omega_c$.
Then our main result in this article
is the following theorem (=Theorem \ref{thm1.1}).
\begin{thm}\label{thm3.1} Let $h\ge 2g$ be a positive integer and $r\ge 2$. Then, for a generic curve $C$ of genus $g$ over an algebraically closed field of characteristic $p>3r$ and
a set $I=\{y_1,\ldots,y_h\}\cup\{z_1,\ldots,z_{2+h-2g}\}$ of generic points on $C$ with $\vec n(y_i)=(r-1,1)$, $\vec n(z_i)=(1,\ldots,1)$, the moduli space $\sU_{C,\,\omega_c}$
is split by a $(p-1)$-power.
\end{thm}

The idea to prove Theorem \ref{thm3.1} is to do induction on genus of the curves and the number of parabolic points. Assume Theorem \ref{thm3.1} holds for $\sU_{\wt X,\,\wt{\omega}_c}$ where $\wt{\omega}_c=(2r,\{\vec n(x), \vec a_c(x)\}_{I\cup \{x_1,x_2\}})$. Let $X$ be the projective curve with one node $x_0$ by gluing $x_1$ and $x_2$. Then there is a moduli space $\sU_{X,\,\omega_c}$ of semistable parabolic torsion free sheaves on $X$. Let $\sB_X\subset \sU_{X,\,\omega_c}$ be the open set of parabolic bundles and $\sB^{L_0}_X\subset \sB_X$ be the
closed subset consisting of parabolic bundles with fixed determinant $L_0\in J^d_X$. Then main results of next section will imply

\begin{prop}\label{prop3.2} If Theorem \ref{thm3.1} holds for $\sU_{\wt X,\,\wt{\omega}_c}$. Then
\begin{itemize} \item [(1)] ${\rm H}^1(\sU_{X,\,\omega_c},\Theta_{\sU_{X,\,\omega_c}})=0$;
\item [(2)] There exists a section $\sigma\in {\rm H}^0(\sU_{X,\,\omega_c},\Theta_{\sU_{X,\,\omega_c}})$ such that $\sigma^{p-1}|_{\sB_X^{L_0}}$
defines a Frobenius splitting of a nonempty open subset  of $\sB_X^{L_0}$ for generic $L_0\in J^d_X$.
\end{itemize}
\end{prop}

Assume Proposition \ref{prop3.2} holds, we prove Theorem \ref{thm3.1} by degeneration arguments. Take a flat family $C_S\rightarrow S$ of (genus $g$) projective curves over a smooth curve $S$ such that a fiber $C_{s_0}=X$ and $C_s$ ($s\neq s_0$) are smooth projective curves of genus $g$. Fix a line bundle $\sO(1)$ over $C_S$ such that
$\mathrm{deg}(\sO(1)|_{C_s})=c>0$ for all $s\in S$, let $\chi=d+r(1-g)$, $P$ denote the polynomial
$P(m)=crm+\chi$ and $V=k^{\oplus P(N)}$. Assume that $N$ is large enough such that $F(N)$ is generated by global sections and ${\rm H}^1(F(N))=0$ for any semistable parabolic sheaf $F$ of degree $d$ on $C_s$ ($s\in S$).
Let $\bold{Q}_S$ be the relative Quot scheme of quotients $V\otimes \sO_{C_S}(-N)\rightarrow F\rightarrow 0$ of rank $r$ and
degree $d$ on $C_S$ and
$$V\otimes \sO_{C_S\times \bold{Q}_S}(-N)\rightarrow \sF\rightarrow 0$$
be universal quotient on $C_S\times \bold{Q}_S$. Given sections $S_x:S\rightarrow C_S$ ($x\in I$) with $S_x\cap S_y=\emptyset$ for $x\neq y\in I$, let $\mathcal{F}_{S_x}=\mathcal{F}_S|_{S_x\times \bold{Q}}$ and ${\rm Flag}_{\vec{n}(x)}(\mathcal{F}_{S_x})$ be the relative flag scheme of type $\vec n(x)$. Let
$$\mathcal{R}_{S}=\underset{x\in I}{\times}{_{\bold{Q}_{S}} {\rm Flag}_{\vec{n}(x)}(\mathcal{F}_{S_x})},$$ and $\mathcal{R}^{ss}_{S,\,\omega}$ (resp. $\mathcal{R}^s_{S,\,\omega}$) be the open subset of GIT semi-stable (resp. GIT stable) points (respect to $\omega_c$).
Then there exist GIT quotients
$$\mathcal{U}_{C_S,\,\omega_c}:=\mathcal{R}^{ss}_{S,\,\omega_c}//\mathrm{SL(V)}\,\,\,\,\,\,(\text{resp.}\, \mathcal{U}_{C_S,\,\omega_c}^s:=\mathcal{R}^{s}_{S,\,\omega_c}//\mathrm{SL(V)}).$$

\begin{prop}\label{prop3.3} There exists a
$S$-flat scheme $\sU_{C_S,\,\omega_c}\to S$ such that scheme-theoretic fiber at $s\in S$ is isomorphic to moduli space $\mathcal{U}_{C_s,\, \omega_c}$ of
semi-stable parabolic sheaves $E$ of rank $r$ and degree $d$ on $C_s$ with parabolic structures at points $\{s_x\}_{x\in I}$ determined by $\omega_c$. Moreover, there is a relative
ample line bundle $\Theta_{\sU_{C_S,\,\omega_c}}$ on $\sU_{C_S,\,\omega_c}$ such that
$$\Theta_{\sU_{C_S,\,\omega_c}}|_{\mathcal{U}_{C_s,\, \omega_c}}=\Theta_{\mathcal{U}_{C_s,\, \omega_c}}\,\,\,(\forall\,\,s\in S).$$
\end{prop}

\begin{proof} By the choice of $\omega_c$ and Lemma \ref{lem2.13}, the descent lemma holds and the GIT quotient $\mathcal{U}_{C_S,\,\omega}$ commutes
with base changes.
\end{proof}

\begin{cor}\label{cor3.4} Assume that Proposition \ref{prop3.2} holds. Then Theorem \ref{thm3.1} holds if it holds for $g=0$.
\end{cor}

\begin{proof} Let $\sB_{C_S}\subset \sU_{C_S,\,\omega_c}$ be the open set consisting of locally free sheaves and $\mathrm{Det}: \sB_{C_S}\rightarrow J_{C_S}^d$ be determinant map to the relative Jacobian variety. Let $\sL$ be a line bundle on $C_S$ of degree $d$ on each fiber $C_s$ ($\forall\,s\in S$). Then we have a flat family
$\sB_{C_S}^{\sL}\to S$ with fiber $\sB_{C_s}^{\sL_s}=\sU_{C_s,\,\omega_c}^{\sL_s}$ for $s\neq s_0$ and $\sB_{C_{s_0}}^{\sL_{s_0}}:=\sB_X^{L_0}\subset \sU_{X,\,\omega_c}$ is the locally closed
subvariety of parabolic bundles with fixed determinant $\sL_{s_0}=L_0$ on $X$.

Since ${\rm H}^1(\sU_{X,\,\omega_c},\Theta_{\sU_{X,\,\omega_c}})=0$, the section $\sigma\in {\rm H}^0(\sU_{X,\,\omega_c},\Theta_{\sU_{X,\,\omega_c}})$ in Proposition \ref{prop3.2} spreads to
$\mathcal{U}_{C_s,\, \omega_c}$ for $s$ in an open neighbourhood of $s_0\in S$. Replace $S$ by an open neighbourhood of $s_0\in S$, we have
$$\wt{\sigma}\in {\rm H}^0(\sU_{C_S,\,\omega_c},  \Theta_{\sU_{C_S,\,\omega_c}})$$
such that $\wt{\sigma}|_{\mathcal{U}_{C_{s_0},\, \omega_c}}=\sigma\in {\rm H}^0(\sU_{X,\,\omega_c},\Theta_{\sU_{X,\,\omega_c}})$. Then $\tau=\wt{\sigma}|_{\sB_{C_S}^{\sL}}$ is a section
on the family $\sB_{C_S}^{\sL}\to S$ such that
$$\tau|_{\sB_{C_s}^{\sL_s}}\in {\rm H}^0(\sU_{C_s,\,\omega_c}^{\sL_s},\, \omega^{-1}_{\sU_{C_s,\,\omega_c}^{\sL_s}})\,\,(s\neq s_0),
\qquad \tau|_{\sB_{C_{s_0}}^{\sL_{s_0}}}=\sigma|_{\sB^{L_0}_X}.$$

Now, by (2) of Proposition \ref{prop3.2}, the section $\tau^{p-1}$ splits a nonempty open set of
the fiber $\sB_{C_{s_0}}^{\sL_{s_0}}$. Then, by (4) of Proposition \ref{prop2.2}, $\tau^{p-1}$ splits the fibers
$\sB_{C_s}^{\sL_s}=\sU_{C_s,\,\omega_c}^{\sL_s}$ for $s\neq s_0$ in a neighbourhood of $s_0\in S$. Thus, by using Lemma 5.4 of
\cite{M.R}, we are done.
\end{proof}

To show how Proposition \ref{prop3.2} follows main results of next section, we have to recall the normalization of moduli space $\sU_{X,\,\omega}$
of semistable parabolic torsion free sheaves of rank $r$ and degree $d$ on $X$. Let $\wt X$ be an irreducible smooth projective curve of genus $g-1$ and $I\cup\{x_1,\,x_2\}$ be a subset of points of $\wt X$. By identifying $x_1$ and $x_2$, we get an
irreducible projective curve $X$ of (arithmetic) genus $g$ with one node $x_0\in X$. Let
\ga{3.1}{\pi:\wt X\to X,\quad \pi^{-1}(x_0)=\{x_1,\,x_2\} } be the normalization morphism, we identify $I\subset \wt X$ with a subset $I\subset X$ of smooth points.
Let $\omega=(k,\{\vec n(x), \vec a(x)\}_{x\in I})$ be any given data and $\sU_{X,\,\omega}$ (resp. $\sU_{\wt X,\,\omega}$) be the moduli space
of semistable parabolic torsion free sheaves of rank $r$ and degree $d$ on $X$ (resp. $\wt X$) with parabolic structures determined by $\omega$. $\sU_{\wt X,\,\omega}$ and $\sU_{X,\,\omega}$ are related through a moduli space $\sP_{\omega}=\{\,(E,Q)\,\}$ of semistable generalized parabolic sheaves (GPS) on $\wt X$,
where $E$ is a parabolic sheaf of rank $r$ and degree $d$ on $\wt X$ with parabolic structures determined by $\omega$
and $Q$ is a $r$-dimensional quotient $E_{x_1}\oplus E_{x_2}\xrightarrow{q} Q$, whose construction we recall below.

\begin{defn}\label{defn3.5} Let $0<\alpha\le 1$ be a rational number.
 A GPS $(E,Q)$ on $\wt X$ is called $\alpha$-\emph{semistable} (resp.
$\alpha$-\emph{stable}), if for nontrivial subsheaf $F\subset E$ such that
$E/F$ is torsion free outside $\{x_1,x_2\},$ we have
$$\mathrm{par}\chi(F)-\mathrm{dim}(Q^{F})\cdot \alpha\leq
r(F)\cdot\frac{\mathrm{par}\chi(E)-\mathrm{dim}(Q)\cdot \alpha}{r(E)} \,\quad (\text{resp.
$<$}),$$ where $Q^{F}=q(F_{x_1}\oplus F_{x_2})\subset Q.$
\end{defn}

\begin{thm}\label{thm3.6} For any $\omega=(k, \{\vec n(x), \vec a(x)\}_{x\in I})$ and $\alpha\in (0,1]$, there exists
a normal projective variety $\sP_{\omega,\,\alpha}$, which is the
coarse  moduli space of $s$-equivalence classes of $\alpha$-semistable GPS on $\wt X$ with parabolic structures given by $\omega$.  The set of $\alpha$-stable GPS forms an open subset $$\mathcal{P}^s_{\omega,\,\alpha}\subset \mathcal{P}_{\omega,\,\alpha}$$ which is a smooth variety.
\end{thm}

\begin{proof} Let $\wt{\bold Q}$ be the Quot scheme of quotients $\wt V\otimes\sO_{\wt X}(-N)\to \wt F\to 0$ (of rank
$r$ and degree $d$) on $\wt X$ and $\wt V\otimes\sO_{\wt X\times\wt{\bold Q}}(-N)\to \wt\sF\to 0$ be the universal quotient on $\wt X\times\wt{\bold Q}$.
Let $\wt\sF_x=\wt\sF|_{\{x\}\times\wt{\bold Q}}$ and
$$\wt\sR=\underset{x\in I}{\times_{\wt{\bold Q}}}{\rm Flag}_{\vec n(x)}(\wt\sF_x)\to \wt{\bold Q},\quad \wt\sR':={\rm Grass}_r(\wt\sE_{x_1}\oplus\wt\sE_{x_2})\xrightarrow{\rho} \wt\sR$$
where ${\rm Flag}_{\vec n(x)}(\wt\sF_x)\to\wt{\bold Q}$ is the relative flag scheme of type $\vec n(x)$ and $\wt\sE$ is the pullback of $\wt\sF$.
Let $ \wt\sE_{x_1}\oplus\wt\sE_{x_2}\to \sQ\to 0$ be the universal quotient where $\wt\sE_{x_i}=\wt\sE|_{\{x_i\}\times\wt\sR}$.
The construction of $\sP_{\omega,\,\alpha}$ is already given in \cite{B1} for $0<\alpha<1$ (see also \cite{Su3} for $\alpha=1$). In fact, there is a $\mathbb{Q}$-polarization $\hat{\Theta}_{\omega,\alpha}$ on $\wt\sR'$
such that the open set ${\wt\sR}^{'ss}_{\omega,\alpha}\subset \wt\sR'$ (resp. ${\wt\sR}^{'s}_{\omega,\alpha}\subset \wt\sR'$) of GIT semistable points (resp. GIT stable points)
respect to the polarization $\hat{\Theta}_{\omega,\alpha}$ is precisely the open set of $\alpha$-\emph{semistable} (resp.
$\alpha$-\emph{stable}) GPS on $\wt X$. Then
$\sP_{\omega,\,\alpha}$ and $\sP^s_{\omega,\,\alpha}$ are the GIT quotients
\ga{3.2}{{\wt\sR}^{'ss}_{\omega,\alpha}\to{\wt\sR}^{'ss}_{\omega,\alpha}//{\rm SL}(\wt V),\quad {\wt\sR}^{'s}_{\omega,\alpha}\to {\wt\sR}^{'s}_{\omega,\alpha}//{\rm SL}(\wt V).}
\end{proof}

When $0<\alpha<1$, ${\wt\sR}^{'ss}_{\omega,\alpha}\subset {\wt\sR}^{'}_F$ where ${\wt\sR}^{'}_F\subset {\wt\sR}^{'}$ is the open set of GPS $(E,Q)$ with $E$ being locally free. When $\alpha=1$, for $(E,Q)\in {\wt\sR}^{'ss}_{\omega,1}:={\wt\sR}^{'ss}_{\omega}$, $E$ may have torsion at $x_1$ and $x_2$. But there is an irreducible open set
$\sH\subset {\wt\sR}^{'}$ (which is normal) such that ${\wt\sR}^{'ss}_{\omega}\subset \sH$ for any $\omega$ (see Notation 2.1 and Remark 2.1 of \cite{S1}). In fact, it is easy to see ${\wt\sR}^{'ss}_{\omega,\alpha}\subset\sH$ for any $0< \alpha\le 1$. Moreover, by Proposition 5.3 of \cite{Su3} (see \cite[Proposition 5.2]{S1} for detail), we have
\ga{3.3}{{\rm Codim}(\sH\setminus {\wt\sR}^{'ss}_{\omega})>\underset{0<r_1<r}{\mathrm{min}}\{r_1(r-r_1)(g-1)+\Sigma_{x\in I}\Sigma_{x,r_1}(\omega)\}.}
For any $(E,Q)\in {\wt\sR}^{'ss}_{\omega}$,  $\phi((E,Q)):=\mathrm{ker}\{ \pi_*E\to\, _{x_0}(\pi_*E)_{x_0}\to\,_{x_0}Q\}$ is
a semistable parabolic torsion free sheaf of rank $r$ and degree $d$ on $X$ with parabolic structures determined by $\omega$, which defines the normalization morphism
(see \cite[Proposition 2.1]{S1} for  detail)
\ga{3.4} {\phi: \sP_{\omega}:=\sP_{\omega,\,1}\to\sU_{X,\,\omega}.}
For $(E,Q)\in {\wt\sR}^{'}$, let $E_{x_i}\xrightarrow{q_i}Q$ ($i=1,\,2$) be the homomorphism induced by the quotient homomorphism $E_{x_1}\oplus E_{x_2}\xrightarrow{q}Q\to 0$,
\ga{3.5} {\hat{D}_{F,i}=\{\,(E,Q)\,|\,\text{$E_{x_i}\xrightarrow{q_i} Q$ is not an isomorphism}\,\}}
and $\hat{D}_i\subset \sH$ be the Zariski closure of $\hat{D}_{F,i}$. Then they induce effective divisors $D_i\subset\sP_{\omega}$. We collect some facts from \cite{N.R} and \cite{S1} where their proofs work in our case (see Remark \ref{rmk2.14}).
\begin{prop}\label{prop3.7} Let $\sW\subset \sU_{X,\,\omega}$ be the closed subset of non-locally free sheaves and $\sI_{\sW}\subset \sO_{\sU_{X,\,\omega}}$ be its idea sheaf. Then
\begin{itemize}
\item [(1)] $\sW\subset \sU_{X,\,\omega}$ is the non-normal locus of $\sU_{X,\,\omega}$ and $$\phi^{-1}(\sW)=D_1+D_2.$$
\item [(2)] $\sU_{X,\,\omega}$ is seminormal and $\sI_{\sW}=\phi_*\sO_{\sP_{\omega}}(-D_1-D_2)$. In particular, for any line bundle $\sL$ on $\sU_{X,\,\omega}$,
$${\rm H}^0(\sP_{\omega}, \phi^*\sL(-D_1-D_2))={\rm H}^0(\sU_{X,\,\omega}, \sI_{\sW}\otimes\sL).$$
\item [(3)] When $p>3r$, the homomorphism $$\phi^*:{\rm H}^1(\sU_{X,\,\omega_c},\Theta_{\sU_{X,\,\omega_c}})\hookrightarrow {\rm H}^1(\sP_{\omega},\Theta_{\sP_{\omega_c}})$$
is injective where $\omega_c$ is canonical weight and
$\Theta_{\sP_{\omega_c}}=\phi^*\Theta_{\sU_{X,\,\omega_c}}.$
\end{itemize}
\end{prop}

\begin{proof} (1) This was proved in Proposition 2.1 of \cite{S1}. (2) Seminormality of $\sU_{X,\,\omega}$ for rank two was proved in \cite{N.R} (Theorem 3 of \cite{N.R}) and the case of arbitrary rank was proved in \cite{S1} (Theorem 4.2 of \cite{S1}), the equality of $\sI_{\sW}=\phi_*\sO_{\sP_{\omega}}(-D_1-D_2)$ follows the general property of seminormal varieties (see \cite[Lemma 3.7]{N.R}). (3) This is Lemma 5.5 of \cite{S1} since the exact sequence (4.2) in the proof of Proposition 4.2 of \cite{S1} still
splits when $p>3r$ (see Lemma \ref{lem2.16}). In fact, let ${\widetilde{\sR}}^s_{\widetilde{X},\,\omega_c}\subset \wt \sR_F$ be the open subset consisting of stable parabolic bundles respect to $\omega_c$, and $\rho:\rho^{-1}({\widetilde{\sR}}^s_{\widetilde{X},\,\omega_c})\rightarrow {\widetilde{\sR}}^s_{\widetilde{X},\,\omega_c}$. The Grassmannian $\mathrm{Grass}_r(\sF_{x_1}\oplus \sF_{x_2})$ (over $\mathcal{U}^s_{\widetilde{X},\,\omega_c}$) is a principal $\mathrm{PGL}(\widetilde{V})$-bundle over $\rho^{-1}({\widetilde{\sR}}^s_{\widetilde{X},\,\omega_c})$. Then the sequence (4.2) in the proof of Proposition 4.2 of \cite{S1} can be replaced by the similar sequence over $\mathcal{U}^s_{\widetilde{X},\,\omega_c}$.
\end{proof}

The following theorem will be proved in next section.

\begin{thm}\label{thm3.8} For any $\omega_c=(2r,\{\vec n(x), \vec a_c(x)\}_{x\in I})$ such that $\sU^s_{\wt X,\,\omega_c}$ is
nonempty, let $\wt\omega_c$ be the canonical weight
determined by $\{\vec n(x)\}_{x\in I\cup\{x_1,\,x_2\}}$ where $\vec n(x_i)=(1,\ldots, 1)$ $(i=1,\,2)$.
Assume that $\sU_{\wt X,\,\wt\omega_c}$ is split by a $(p-1)$-power and $p>3r$. Then there exists a section
$$\wt\sigma\in {\rm H}^0(\sP_{\omega_c}, \omega^{-1}_{\sP_{\omega_c}})$$ such that
$\wt\sigma$ vanishes on $D_1+D_2$ and $\wt\sigma^{p-1}$ splits $\sP_{\omega_c}$.
\end{thm}

\begin{proof}[Proof of Proposition \ref{prop3.2}] By Theorem \ref{thm3.8}, ${\rm H}^1(\sP_{\omega_c}, \Theta_{\sP_{\omega_c}})=0$, which implies
${\rm H}^1(\sU_{X,\,\omega_c},\Theta_{\sU_{X,\,\omega_c}})=0$ by Proposition \ref{prop3.7} (3) when $p>3r$. This proves (1) of Proposition \ref{prop3.2}.
To prove (2) of Proposition \ref{prop3.2}, let ${\rm Det}:\sP_{\omega_c}\to J^d_{\wt X}$ be the determinant morphism, we have
\ga{3.6}{ \omega^{-1}_{\sP_{\omega_c}}=\Theta_{\sP_{\omega_c}}\otimes {\rm Det}^*(\Theta_{J^d_{\wt X}}^{-1}).}
This can be deduced as follows. Let $\omega$ be a minor modification of $\omega_c$ such that moduli space $\sP_{\omega,\,1^-}$ is fine
and $\widetilde{\sR}^{'s}_{\omega,\,1^-}\rightarrow \sP_{\omega,\,1^-}$ is a principal $\mathrm{PGL}(\widetilde{V})$-bundle. Then, by Proposition 4.6 of \cite{SZ}, we have \ga{3.7}{\omega^{-1}_{\sP_{\omega,\,1^-}}=f^{\ast}(\Theta_{\sP_{\omega_c}})\otimes {\rm Det}^*(\Theta_{J^d_{\wt X}}^{-1})} where $f:\sP_{\omega,\,1^-}\rightarrow \sP_{\omega_c}$ is a birational morphism (see Proposition 6.3 (2) of \cite{S4}). Let $U\subset \sP_{\omega_c}$ be the largest open set such that $f^{-1}$ is well-defined. Then \eqref{3.6} holds on $U$ by restricting \eqref{3.7} to $f^{-1}(U)\xrightarrow{\backsim} U$, which implies that \eqref{3.6} holds on $\sP_{\omega_c}$ since $\sP_{\omega_c}\setminus U$ has codimension at least two.
By tensoring  a section
of ${\rm Det}^*(\Theta_{J^d_{\wt X}})$, we get a $${\wt\sigma}'\in {\rm H}^0(\sP_{\omega_c},\Theta_{\sP_{\omega_c}}(-D_1-D_2))$$
such that $({\wt\sigma}')^{p-1}|_{\sP^{\wt L}_{\omega_c}}$ splits $\sP^{\wt L}_{\omega_c}:={\rm Det}^{-1}(\wt L)$ for generic $\wt L\in J^d_{\wt X}$ (see Proposition \ref{prop2.2} (2)). Thus, by (2) of Proposition \ref{prop3.7}, there is a
$\sigma\in {\rm H}^0(\sU_{X,\,\omega},\Theta_{\sU_{\wt X,\omega_c}})$ such that ${\wt\sigma}'=\phi^*(\sigma)$.

Let $L\in J^d_X$, $\wt L=\pi^*(L)$, $\sP^{\wt L}_B:=\sP^{\wt L}_{\omega_c}\setminus\{D_1\cup D_2\}$ and ${\rm Isom}(\wt L_{x_1},\wt L_{x_2})$ be the open set of ${\rm Hom}(\wt L_{x_1},\wt L_{x_2})$
consisting of isomorphisms. Then there is a morphism $\Delta:\sP^{\wt L}_B\to {\rm Isom}(\wt\sL_{x_1},\wt\sL_{x_2})$, its generic fiber $\Delta^{-1}(\eta)\subset \sP^{\wt L}_B$ has a nonempty open subset which is split by $({\wt\sigma}')^{p-1}|_{\Delta^{-1}(\eta)}$ (see Proposition \ref{prop2.2} (5)). The point $\eta\in {\rm Isom}(\wt L_{x_1},\wt L_{x_2})$ determines a line bundle $L_0\in J^d_X$ such that
$\pi^*(L_0)=\wt L=\pi^*(L)$. Then the isomorphism $\sP_B:=\sP_{\omega_c}\setminus\{D_1\cup D_2\}\xrightarrow{\phi}\sB_X=\sU_{X,\omega_c}\setminus\sW$ induces an isomorphism
$$\Delta^{-1}(\eta)\xrightarrow{\phi}\sB_X^{L_0},$$
$\sigma^{p-1}|_{\sB_X^{L_0}}$ splits a nonempty open subset of  $\sB_X^{L_0}$ since $({\wt\sigma}')^{p-1}|_{\Delta^{-1}(\eta)}=\phi^*(\sigma^{p-1}|_{\sB_X^{L_0}})$ splits a nonempty open subset of $\Delta^{-1}(\eta)$.
\end{proof}

Thus, under assumption that Theorem \ref{thm3.8} holds, the main theorem is reduced to the proof of $g=0$ case (=Theorem \ref{thm1.3}).
\begin{thm}\label{thm3.9} For any $h>0$, let $I\subset \mathbb{P}^1$ be a set of $2+2h$ points with
$\vec n(y_1)=\cdots=\vec n(y_{h})=(r-1,1)$ $(y_i\in I)$ and  $\vec n(x)=(1,\ldots, 1)$ $(\forall x\in I-\{y_1,\ldots, y_h\})$.
Then, when $p>3r$, the moduli space $\mathcal{U}_{\mathbb{P}^1,\,{\omega}_c}$ is split by a $(p-1)$-power for generic choice of points in $I$.
\end{thm}

We prove Theorem \ref{thm3.9} by induction on $h$. The case $h=1$ will be proved in Section 5. Assume $h>1$ and Theorem \ref{thm3.9}
holds for $h-1$. Let $X:=X_1\cup X_2$ where $X_i=\mathbb{P}^1$ meeting at only one nodal point $x_0$ and
$\pi:\widetilde{X}=X_1\sqcup X_2\rightarrow X$ be its normalization with $\pi^{-1}(x_0)=\{x_1,x_2\}$. Let $I_1=\{y_h,z_1,z_2\}\subset X_1$ and
$I_2=\{y_1,...,y_{h-1},z_3,...,z_{h+2}\}\subset X_2$ with $\vec n(y_i)=(r-1,1)$ and $\vec n(z_i)=(1,...,1)$. Then there is a moduli space $\mathcal{U}_{X,\,\omega_c}$ of semistable parabolic torsion free sheaves on $X$ with parabolic structures given by $\omega_c=(2r,\{\vec n(x),\vec a_c(x)\}_{x\in I_1\cup I_2}\})$ (see \cite{S4} or \cite{Su3} for the construction).
Let $\sB_{X}\subset \mathcal{U}_{X,\,\omega_c}$ be the open subset of parabolic bundles, then we have
\begin{prop}\label{prop3.10} Let $\omega_{c_i}=(2r,\{\vec n(x),\vec a_c(x)\}_{x\in I_i\cup\{x_i\}})$ be the canonical weight determined by $\{\vec n(x)\}_{x\in I_i\cup\{x_i\}}$ where $\vec n(x_i)=(1,...,1)$. If Theorem \ref{thm3.9} holds for $\sU_{X_i,\,\omega_{c_i}}$ $(i=1,2)$. Then
\begin{itemize} \item [(1)] ${\rm H}^1(\sU_{X,\,\omega_c},\Theta_{\sU_{X,\,\omega_c}})=0$;
\item [(2)] There exists a section $\sigma\in {\rm H}^0(\sU_{X,\,\omega_c},\Theta_{\sU_{X,\,\omega_c}})$ such that $\sigma^{p-1}|_{\sB_X}$
defines a Frobenius splitting of $\sB_X$.
\end{itemize}\end{prop}

Assume Proposition \ref{prop3.10} holds, we will prove Theorem \ref{thm3.9} by degeneration arguments. Take a flat family $C_S\rightarrow S$ of curves over a smooth $S$ with fibers $C_{s_0}=X$, $C_s=\mathbb{P}^1$ ($s\neq s_0$) and disjoint sections $S_x:S\rightarrow C_S,\,\, (x\in I_1\cup I_2)$ such that $I_i=\{S_x\cap X_i\}_{x\in I_i}$. Then we have

\begin{prop}\label{prop3.11}
There is a $S$-flat scheme $\mathcal{U}_{C_S,\,\omega_c}\rightarrow S$ such that scheme-theoretic fibers at closed points $s\in S$ are isomorphic to moduli spaces $\mathcal{U}_{C_s,\,\omega_c}$ of semistable parabolic sheaves of rank $r$ and degree $d$ on $C_s$ where
$\omega_c=(2r,\{\vec n(S_x\cap C_s),\vec a_c(S_x\cap C_s)\}_{x\in I_1\cup I_2})$ is the canonical weight determined by $\vec n(S_x\cap C_s)=\vec n(x)$ $(x\in I_1\cup I_2)$. Moreover, there is a relative ample line bundle $\Theta_{\mathcal{U}_{C_S},\,\omega_c}$ on $\mathcal{U}_{C_S,\,\omega_c}$ such that $$\Theta_{\mathcal{U}_{C_S},\,\omega_c}|_{\mathcal{U}_{C_s,\,\omega_c}}=\Theta_{\mathcal{
U}_{C_s},\,\omega_c}\,\, (\forall s\in S).$$
\end{prop}
\begin{proof} By the choice of $\omega_c$ and Lemma \ref{lem2.13}, the descent lemma holds and the GIT quotient $\mathcal{U}_{C_S,\,\omega}$ commutes
with base changes.
\end{proof}

\begin{cor}\label{cor3.12} If Proposition \ref{prop3.10} holds and Theorem \ref{thm3.9} holds when $h=1$, then Theorem \ref{thm3.9} holds for any
$h>1$.
\end{cor}
\begin{proof}
This can be done by following a similar argument in the proof of Corollary \ref{cor3.4}, but with a simpler detail. Indeed, the moduli space $\sB_{C_{s_0}}$ naturally has $\omega^{-1}_{\sB_{C_{s_0}}}=\Theta_{\mathcal{U}_{C_S},\,\omega_c}|_{\sB_{C_{s_0}}}$ where $\sB_{C_{s_0}}\subset \mathcal{U}_{C_{s_0},\,\omega_c}$ is the open subset of parabolic bundles.
\end{proof}

To prove Proposition \ref{prop3.10}, we have to recall the normalization of moduli space $\mathcal{U}_{X,\,\omega_c}$.
We fix an ample line bundle $\sO(1)$ of degree $c$ on $\widetilde{X}=X_1\sqcup X_2$ such that $\mathrm{deg}(\sO(1)|_{X_i})=c_i>0$ ($i=1,2$). For any coherent sheaf $E$,
$P(E,n):=\chi(E(n))$ denotes its Hilbert polynomial, which has degree
$1$. We define the rank of $E$ to be
$$r(E):=\frac{1}{\mathrm{deg}(\sO(1))}\cdot \lim \limits_{n\to\infty}\frac{P(E,n)}
{n}.$$ Let $r_i$ denote the rank of the restriction of $E$ to $X_i$
($i=1,2$), then
$$P(E,n)=(c_1r_1+c_2r_2)n+\chi(E),\quad r(E)=
\frac{c_1}{c_1+c_2}r_1+\frac{c_2}{c_1+c_2}r_2.$$ We say that $E$ is
of rank $r$ on $X$ if $r_1=r_2=r$, otherwise it will be said of rank
$(r_1,r_2)$.

Fix a finite set $I=I_1\cup I_2$ of smooth points on $\widetilde{X}$, where
$I_i=\{x\in I\,|\,x\in X_i\}$ ($i=1,2$) and parabolic data $\omega=\{k,\vec
n(x),\vec a(x)\}_{x\in I}$ with
$$\ell:=\frac{k\chi-\sum_{x\in I}\sum^{l_x}_{i=1}d_i(x)r_i(x)}{r}$$
(recall $d_i(x)=a_{i+1}(x)-a_i(x)$, $r_i(x)=n_1(x)+\cdots+n_i(x)$). Let
\ga{3.8}{n^{\omega}_j=\frac{1}{k}\left(r\frac{c_j}{c_1+c_2}\ell+\sum_{x\in
I_j}\sum^{l_x}_{i=1}d_i(x)r_i(x)\right)\,\,\,(j=1,\,\,2).}

\begin{defn}\label{defn3.13} (1) For any coherent sheaf $F$ of rank $(r_1,r_2)$, let
$$m(F):= \frac{r(F)-r_1}{k}\sum_{x\in I_1}a_{l_x+1}(x)+
\frac{r(F)-r_2}{k}\sum_{x\in I_2}a_{l_x+1}(x),$$ the modified
parabolic Euler characteristic and slope of $F$ are
$${\rm par}\chi_m(F):={\rm par}\chi(F)+m(F),\quad {\rm par}\mu_m(F):=\frac{{\rm par}\chi_m(F)}{r(F)}.$$

(2) Let $0<\alpha \leq 1$ be a rational number.  A GPS $(E,E_{x_1}\oplus E_{x_2}\xrightarrow{q}Q)$ is called $\alpha$-semistable (resp. stable), if for every nontrivial subsheaf $E'\subset E$ such that
$E/E'$ is torsion free outside $\{x_1,x_2\},$ one has, with the
induced parabolic structures at points $\{x\}_{x\in I}$,
$$\mathrm{par}\chi_m(E')-\mathrm{dim}(Q^{E'})\cdot \alpha\leq
r(E')\cdot\frac{\mathrm{par}\chi_m(E)-\mathrm{dim}(Q)\cdot \alpha}{r(E)} \,\quad (\text{resp.
$<$}),$$ where $Q^{E'}=q(E'_{x_1}\oplus E'_{x_2})\subset Q.$
\end{defn}
\begin{thm}[see Theorem 2.1 of \cite{S4} or Theorem 2.26 of \cite{Su3}]\label{thm3.14} For any data $\omega=(\{k,\vec n(x),\,\,\vec a(x)\}_{x\in I_1\cup I_2},\sO(1))$ and $\alpha\in (0,1]$,
the coarse  moduli space $\sP_{\omega,\,\alpha}$ of $s$-equivalence classes of $\alpha$-semistable GPS on $\wt X$ with parabolic structures
at points $\{x\}_{x\in I}$ given by the data $\omega$ is a disjoint union of at most $r+1$ irreducible, normal projective varieties $\sP_{\chi_1,\chi_2}$
$(\chi_1+\chi_2=\chi+r,\, n_j^{\omega}\le\chi_j\le n_j^{\omega}+r)$. The set of $\alpha$-stable GPS forms an open subset $$\sP^s_{\omega,\,\alpha}\subset \sP_{\omega,\,\alpha}$$ which is a smooth variety.
\end{thm}
\begin{proof} The case $\alpha=1$ is already given in \cite{S4}, and general case can be proved by following a similar argument. When $0<\alpha <1$, every semistable GPS $(E,Q)$ should also satisfy
$$n_1^{\omega}\leq \chi(E_1)\leq n_1^{\omega}+r,\,\,n_2^{\omega}\leq \chi(E_2)\leq n_2^{\omega}+r$$ (see (2.2) of \cite{S4}) where $E_1=E|_{X_1}$ and $E_2=E|_{X_2}$.
For fixed $\chi_1$, $\chi_2$ satisfying $\chi_1+\chi_2=\chi+r$ and $n_j^{\omega}\le\chi_j\le n_j^{\omega}+r$ ($j=1,\,2$), we give briefly the construction
of $\sP_{\omega}=\sP_{\chi_1,\chi_2}$. Let
$$P_i(m)=c_irm+\chi_i,\quad \sW_i=\sO_{X_i}(-N),\quad V_i=k^{P_i(N)},$$ where
$\sO_{X_i}(1)=\sO (1)|_{X_i}$ has degree $c_i$. Consider the Quot schemes $\textbf{Q}_i=\mathrm{Quot}(V_i\otimes\sW_i,
P_i)$, the universal quotient
$V_i\otimes\sW_i\to \sF_i\to 0$ on $X_i\times \textbf{Q}_i$ and
the relative flag scheme
$$\sR_i=\underset{x\in I_i}{\times_{\textbf{Q}_i}}
{\rm Flag}_{\vec n(x)}(\sF_{i,x})\to \textbf{Q}_i.$$ Let
$\sF=\sF_1\oplus\sF_2$ denote direct sum of pullbacks of $\sF_1$,
$\sF_2$ on ${\wt X\times
(\textbf{Q}_1\times\textbf{Q}_2),}$
 $\sE$ be the pullback of $\sF$ to $\wt
X\times(\sR_1\times\sR_2)$, and
\ga{3.9}{\widetilde{\sR}={\rm Grass}_r(\sE_{x_1}\oplus\sE_{x_2})\xrightarrow{\rho}\sR=\sR_1\times\sR_2\to
\textbf{Q}=\textbf{Q}_1\times\textbf{Q}_2.}
For given data $\omega=(\{k,\vec n(x),\,\vec a(x)\}_{x\in I_1\cup I_2},\mathcal{O}(1))$ and $\alpha\in (0,1]$, there is a $\mathbb{Q}$-polarization $\hat{\Theta}_{\omega,\,\alpha}$ on $\widetilde{\sR}$
such that the open subset $\widetilde{\sR}^{ss}_{\omega,\,\alpha}\subset \widetilde{\sR}$ (resp. $\widetilde{\sR}^s_{\omega,\, \alpha}\subset \widetilde{\sR}$) of GIT semistable points (resp. GIT stable points) under the action of $G=(\mathrm{GL}(V_1)\times \mathrm{GL}(V_2))\cap \mathrm{SL}(V_1\oplus V_2)$ on $\widetilde{\sR}$ respect to the polarization $\hat{\Theta}_{\omega,\,\alpha}$ is precisely the open set of $\alpha$-semistable (resp. $\alpha$-stable) GPS on $\widetilde{X}$. Then $\sP_{\chi_1,\chi_2}$ and $\sP^s_{\chi_1,\chi_2}$ are the GIT quotients
$${\wt\sR}^{ss}_{\omega,\alpha}\to\sP_{\chi_1,\chi_2}:= {\wt\sR}^{ss}_{\omega,\alpha}//G,\quad {\wt\sR}^{s}_{\omega,\alpha}\to\sP^s_{\chi_1,\chi_2}:= {\wt\sR}^{s}_{\omega,\alpha}//G.$$
\end{proof}
\begin{rmk}\label{rmk3.15} Every $\alpha$-semistable GPS $(E,E_{x_1}\oplus E_{x_2}\xrightarrow{q}Q)$ satisfies $$\aligned &n_1^{\omega}+(r-\mathrm{dim}(Q^{E_2}))\cdot \alpha \leq \chi(E_1)\leq n_1^{\omega}+\mathrm{dim}(Q^{E_1})\cdot \alpha \\& n_2^{\omega}+(r-\mathrm{dim}(Q^{E_1}))\cdot \alpha \leq \chi(E_2)\leq n_2^{\omega}+\mathrm{dim}(Q^{E_2})\cdot \alpha\endaligned$$ where $E_i=E|_{X_i}$. Assume $n_j^{\omega}$ ($j=1,2$) are not integers, moduli spaces $\mathcal{P}_{\omega}$ ($\alpha=1$) and $\mathcal{P}_{\omega,\,1^-}$ both have $r$ irreducible components.
\end{rmk}
When $0\ll \alpha<1$, $\widetilde{\sR}^{ss}_{\omega,\,\alpha}\subset \widetilde{\sR}_F$ where $\widetilde{\sR}_F\subset \widetilde{\sR}$ is the open subset of GPS ($E,Q$) with $E$ being locally free. When $\alpha=1$, for $(E,Q)\in \widetilde{\sR}^{ss}_{\omega,1}:=\widetilde{\sR}^{ss}_{\omega}$, $E$ may have torsion at $x_1$ and $x_2$. But there is an irreducible open subset $\sH\subset \widetilde{\sR}$ (which is normal) such that $\widetilde{\sR}^{ss}_{\omega}\subset \sH$ (see Notation 2.1 of \cite{S4}). For any $(E,Q)\in \widetilde{\sR}^{ss}_{\omega}$, $$\phi((E,Q)):=\mathrm{\mathrm{ker}}\{\pi_{\ast}E\rightarrow \,_{x_0}(\pi_{\ast}E)_{x_0}\rightarrow \,_{x_0}Q\}$$ is a semistable parabolic torsion free sheaf of rank $r$ and degree 0 on $X$ with parabolic structures determined by $\omega$, which defines a morphism $\phi_{\sP_{\chi,\,\chi_2}}:\sP_{\chi_1,\,\chi_2}\rightarrow \mathcal{U}_X$ (see Lemma 3.4 of \cite{S4}).  Running through all the
($\chi_1+\chi_2=\chi+r$, $n_j^{\omega}\le\chi_j\le n_j^{\omega}+r$), there is a morphism ${\phi:=\underset{\chi_1+\chi_2=\chi+r}{\coprod} \phi_{\sP_{\chi_1,\,\chi_2}}:\sP_{\omega}\rightarrow \mathcal{U}_{X,\,\omega}}$
which is the normalization of $\mathcal{U}_{X,\,\omega}$.

Let $\hat{\sD}_i\subset \widetilde{\sR}$ ($i=1,2$)
be the Zariski closure of $\hat{\sD}_{F,i}\subset \widetilde{\sR}_F$ consisting of $(E,Q)$ such that $E_{x_i}\rightarrow Q$ is not an isomorphism, and it induces effective divisor $D_i\subset \sP_{\omega}$. Then we have
\begin{prop}[see \cite{S4} or \cite{SZ1}]\label{prop3.16}
Let $\sW\subset \sU_{X,\,\omega}$ be the closed subset of non-locally free sheaves and $\sI_{\sW}\subset \sO_{\sU_{X,\,\omega}}$ be its idea sheaf. Then
\begin{itemize}
\item [(1)] $\sW\subset \sU_{X,\,\omega}$ is the non-normal locus of $\sU_{X,\,\omega}$ and $$\phi^{-1}(\sW)=D_1+D_2.$$
\item [(2)] $\sU_{X,\,\omega}$ is seminormal and $\sI_{\sW}=\phi_*\sO_{\sP_{\omega}}(-D_1-D_2)$. In particular, for any line bundle $\sL$ on $\sU_{X,\,\omega}$,
$${\rm H}^0(\sP_{\omega}, \phi^*\sL(-D_1-D_2))={\rm H}^0(\sU_{X,\,\omega}, \sI_{\sW}\otimes\sL).$$
\item [(3)] When $p>3r$, the homomorphism $$\phi^*:{\rm H}^1(\sU_{X,\,\omega_c},\Theta_{\sU_{X,\,\omega_c}})\hookrightarrow {\rm H}^1(\sP_{\omega},\Theta_{\sP_{\omega_c}})$$
is injective where $\omega_c$ is canonical weight and
$\Theta_{\sP_{\omega_c}}=\phi^*\Theta_{\sU_{X,\,\omega_c}}.$
\end{itemize}
\end{prop}

In next section, we will prove the following theorem.
\begin{thm}\label{thm3.17} Let $\omega_{c_i}=(2r,\{\vec n(x),\vec a_c(x)\}_{x\in I_i\cup\{x_i\}})$ be canonical weight determined by $\{\vec n(x)\}_{x\in I_i\cup\{x_i\}}$ where $\vec n(x_i)=(1,...,1)$. If Theorem \ref{thm3.9} holds for $\sU_{X_i,\,\omega_{c_i}}$ $(i=1,2)$. Then, when $p>3r$, there is a
$$\tilde{\sigma}\in {\rm H}^0(\sP_{\omega_c},\omega^{-1}_{\sP_{\omega_c}})\quad (\sP_{\omega_c}:=\sP_{\omega_c,1})$$
such that $\tilde{\sigma}$ vanishes on $D_1+D_2$ and $\tilde{\sigma}^{p-1}$ splits $\sP_{\omega_c}$, where
$$\omega_{c}=(2r,\{\vec n(x),\vec a_c(x)\}_{x\in I_1\cup I_2})$$
is the canonical weight determined by $\vec n(x)$ $(x\in\{y_1,...,y_h,z_1,...,z_{h+2}\})$.
\end{thm}

Assume that Theorem \ref{thm3.17} holds, we will prove Proposition \ref{prop3.10}.

\begin{proof}[Proof of Proposition \ref{prop3.10}] By Theorem \ref{thm3.17}, we have $${\rm H}^1(\sP_{\omega_c},\Theta_{\sP,\,\omega_c})=0,$$ which implies
${\rm H}^1(\mathcal{U}_{X,\,\omega_c},\Theta_{\mathcal{U}_X,\,\omega_c})=0$ by Proposition \ref{prop3.16} (3). This is (1) of Proposition \ref{prop3.10}. The
proof of (2) of Proposition \ref{prop3.10} is similar with the proof of Proposition \ref{prop3.2} (2).
\end{proof}

Altogether, to finish the proof of Theorem \ref{thm1.1}, it is enough to show Theorem \ref{thm3.8}, Theorem \ref{thm3.17} and the case $h=1$ of Theorem \ref{thm3.9}.
Theorem \ref{thm3.8} and Theorem \ref{thm3.17} will be proved in Section 4. The case $h=1$ of Theorem \ref{thm3.9} will be proved in Section 5.

\section{$F$-splitting of moduli spaces of generalized parabolic sheaves}

The goal of this section is to prove Theorem \ref{thm3.8} and Theorem \ref{thm3.17}. We prove firstly Theorem \ref{thm3.8} and keep use of notations introduced in last section.
Recall, when $0<\alpha<1$, ${\wt\sR}^{'ss}_{\omega,\alpha}\subset {\wt\sR}^{'}_F$
where ${\wt\sR}^{'}_F\subset {\wt\sR}^{'}$ is the open set of GPS $(E,Q)$ with $E$ being locally free. The divisors
$\hat{D}_{F,i}\subset  {\wt\sR}^{'}_{F}$ ($i=1,\,2$) defined in \eqref{3.5} induce divisors
\ga{4.1} {D_i\subset \sP_{\omega,\alpha}:={\wt\sR}^{'ss}_{\omega,\alpha}//{\rm SL}(\wt V),\quad i=1,\,\,2}
where for simplicity of notation we use $D_i\subset\sP_{\omega,\alpha}$ to denote divisors in different moduli spaces for different $0<\alpha\le 1$.

The following lemma allows us to make minor modifications of $\omega_c$ (see Remark \ref{rmk2.10}). In particular, it is enough to show Theorem \ref{thm3.8} for $\sP_{\omega_c,1^-}$
where $\alpha=1^-<1$ (resp. $\alpha=0^+>0$) is a rational number such that $|1-1^-|$ (resp. $|0^+-0|$) is \textit{sufficiently small}.

\begin{lem}\label{lem4.1} Let $f: X\to Y$ be a proper birational morphism of normal varieties and $D\subset X$ be an effective divisor such that
$f(D)$ be a divisor of $Y$. If there exists $\wt\sigma\in {\rm H}^0(X,\omega^{-1}_{X})$ such that $\wt\sigma|_D=0$ and ${\wt\sigma}^{p-1}$ splits
$X$. Then there exists $\sigma\in {\rm H}^0(Y,\omega^{-1}_{Y})$ such that $\sigma|_{f(D)}=0$ and $\sigma^{p-1}$ splits $Y$.
\end{lem}

\begin{proof} We remark that definitions of $\omega_X^{-1}$ and $F$-splitting of normal varieties are formulated in Remark 1.3.12 of \cite{M.S}.
Let $U\subset Y$ be the largest open set where $f^{-1}:Y\dashrightarrow X$ is well-defined. Then $Y\setminus U$ has codimension at least two since
$Y$ is normal and $f:X\to Y$ is proper. Thus $f: f^{-1}(U)\to U$ is an isomorphism and $D\cap f^{-1}(U)\neq\emptyset$ (since $f(D)$ has codimension one),
which implies that $\wt\sigma|_{f^{-1}(U)}$ induces a section $\sigma\in{\rm H}^0(Y,\omega^{-1}_{Y})$ satisfying the requirements.
\end{proof}

Assume that $|I|$ is large enough and $\sE$ is the universal bundle on $\wt X\times \sU_{\wt X,\,\omega}$ ($\omega$ is generic and a minor modification of $\omega_c$), then $\sP_{\omega,1^-}$ and $\sP_{\omega,0^+}={\rm Grass}_r(\sE_{x_1}\oplus\sE_{x_2})$ are isomorphic in codimension two, and
$$Y:=\mathrm{Flag}_{(1,\ldots, 1)}(\sE_{x_1})\times_{\sU_{\wt X,\,\omega}} \mathrm{Flag}_{(1,\ldots, 1)}(\sE_{x_2})$$
is split by a $(p-1)$-power since $\sU_{\wt X,\,\wt\omega_c}$ is split by a $(p-1)$-power. Then the following construction
implies Theorem \ref{thm3.8} immediately. We will remove the assumption that $|I|$ is large enough in this section. In fact, from the
splitting section $\sigma\in {\rm H}^0(\omega^{-1}_{\sU_{\wt X,\,\wt\omega_c}})$, we will construct a splitting section
$\sigma_0\in {\rm H}^0(\omega^{-1}_{\sP_{\omega,1^-}})$ by the following general construction.

\subsection{A general construction} Let $E_1$, $E_2$ be vector bundles of rank $r$ on a variety $M$,
$\overrightarrow m_j=(\overbrace{1,\ldots, 1}^{j},r-2j,\overbrace{1,\ldots, 1}^{j})$ ($r-j>j$).
Let $\mathrm{Flag}_{\overrightarrow m_j}(E_i)$ ($i=1,2$) be the flag bundle of type $\overrightarrow m_j$ and
$$E_i\rightarrow \mathcal{Q}_{i,r-1}\rightarrow \cdots \mathcal{Q}_{i,r-j}\rightarrow \mathcal{Q}_{i,j}\rightarrow \cdots \rightarrow \mathcal{Q}_{i,1}\rightarrow 0$$ be the universal flag over $M$. For any $k\ge \ell$, define
\ga{4.2} {\mathcal{F}^i_{k,\ell}:={\rm ker}(\mathcal{Q}_{i,k}\twoheadrightarrow \mathcal{Q}_{i,\ell}).}
Let $Gr_0:=\mathrm{Grass}_r(E_1\oplus E_2)\xrightarrow{\rho}M:=\mathrm{Flag}_{\overrightarrow m_0}$,
\ga{4.3}{Gr_j:=\mathrm{Grass}_{r-2j}(\mathcal{F}^{1}_{r-j,j}\oplus \mathcal{F}^{2}_{r-j,j})\xrightarrow{\rho_j} \mathrm{Flag}_{\overrightarrow m_j}}
where $\mathrm{Flag}_{\overrightarrow m_j}:=\mathrm{Flag}_{\overrightarrow m_j}(E_1)\times_M \mathrm{Flag}_{\overrightarrow m_j}(E_2)$ for $1\leq j<\frac{r}{2}$, and
\ga{4.4} {Gr_{\frac{r}{2}}:=\mathrm{Flag}_{(1,\ldots, 1)}(E_1)\times_M \mathrm{Flag}_{(1,\ldots, 1)}(E_2):=\mathrm{Flag}_{\overrightarrow m_{\frac{r}{2}}}}
when $r$ is an even number. Let $\mathcal{F}^{1}_{r-j,j}\oplus \mathcal{F}^{2}_{r-j,j}\xrightarrow{q}\mathcal{Q}_{r-2j}\to 0$
be the universal quotient on $Gr_j$, $D^{[j]}_i\subset Gr_j$ be the divisor where
\ga{4.5} {q_i:\mathcal{F}^{i}_{r-j,j}\subset \mathcal{F}^{1}_{r-j,j}\oplus \mathcal{F}^{2}_{r-j,j}\xrightarrow{q}\mathcal{Q}_{r-2j}} is not an isomorphism
($i=1,\,2$). Then we have

\begin{lem}\label{lem4.2} There exist morphisms $Gr_j\xrightarrow{\varphi_j} Gr_{j-1}$ such that
$$\xymatrix{  &Gr_j\ar[r]^{\tiny \varphi_j}\ar[d]&
     Gr_{j-1}\ar[d] \\
 & \mathrm{Flag}_{\overrightarrow m_j} \ar[r]&
\mathrm{Flag}_{\overrightarrow m_{j-1}}}$$
$(1\le j\le \frac{r}{2})$ are commutative, where $\mathrm{Flag}_{\overrightarrow m_j}\to \mathrm{Flag}_{\overrightarrow m_{j-1}}$ are canonical projections.
Moreover, $\varphi_j$ induce  birational morphisms
$$Gr_j\xrightarrow{\varphi_j} D^{[j-1]}_1\cap D^{[j-1]}_2\subset Gr_{j-1},\quad (1\le j\le\frac{r}{2}),$$
where $D^{[0]}_i:=D_i\subset Gr_0$ $(i=1,\,2)$.
\end{lem}

\begin{proof} When $j=1$, let $\sK_{r-2}:={\rm ker}(\mathcal{F}^{1}_{r-1,1}\oplus \mathcal{F}^{2}_{r-1,1}\xrightarrow{q}\mathcal{Q}_{r-2})$ and
$$\sK_r:={\rm ker}(\sQ^{\vee}_{1,r-1}\oplus \sQ^{\vee}_{2,r-1}\twoheadrightarrow\sF^{1\,\vee}_{r-1,1}\oplus \sF^{2\,\vee}_{r-1,1}\twoheadrightarrow\sK^{\vee}_{r-2}).$$
Then the exact sequence $0\rightarrow \sK_r\rightarrow \sQ^{\vee}_{1,r-1}\oplus \sQ^{\vee}_{2,r-1}\rightarrow \sK^{\vee}_{r-2}\rightarrow 0$ induces
$E_1\oplus E_2\twoheadrightarrow\mathcal{Q}_{1,r-1}\oplus \mathcal{Q}_{2,r-1}\twoheadrightarrow\mathcal{K}_r^{\vee}$, which
defines $$Gr_1\xrightarrow{\varphi_1} Gr_0:=\mathrm{Grass}_r(E_1\oplus E_2).$$
Clearly, $\varphi_1(Gr_1)=D_1\cap D_2$ and $\varphi_1:Gr_1\rightarrow D_1\cap D_2$ is a birational morphism. In fact, $\varphi_1^{-1}$ is well defined on open subset $U_0\subset D_1\cap D_2$ where both $E_1\xrightarrow{q_1} \mathcal{Q}_r$ and $E_2\xrightarrow{q_2} \mathcal{Q}_r$ have rank $r-1$. Indeed, let
$$\mathcal{Q}_{i,r-1}:=q_i(E_i)\subset \sQ_r,\quad \sK_{r-2}:={\rm ker}(\sQ_{1,r-1}\oplus\sQ_{2,r-1}\xrightarrow{a}\sQ_r\to 0),$$
then $\sK_{r-2}\xrightarrow{p_1}\sQ_{1,r-1}$ (induced by projection $\sQ_{1,r-1}\oplus\sQ_{2,r-1}\xrightarrow{p_1}\sQ_{1,r-1}$) is an injective homomorphism
of bundles on $U_0$ since
$$a_2:\sQ_{2,r-1}\subset \sQ_{1,r-1}\oplus\sQ_{2,r-1}\xrightarrow{a}\sQ_r$$
is so on $U_0$. Similarly, $\sK_{r-2}\xrightarrow{p_2}\sQ_{2,r-1}$ is injective. Let
$$\sF^i_{r-1,1}:=p_i(\sK_{r-2})\subset \sQ_{i,r-1},\quad \sQ_{i,1}:=\sQ_{i,r-1}/\sF^i_{r-1,1}\,\,\,(i=1,\,2).$$
Then, on $U_0$, we have flags $E_i\twoheadrightarrow \sQ_{i,r-1}\twoheadrightarrow\sQ_{i,1}$ with
$$\sF^i_{r-1,1}={\rm ker}(\sQ_{i,r-1}\twoheadrightarrow\sQ_{i,1})$$
and $\mathcal{F}^{1}_{r-1,1}\oplus \mathcal{F}^{2}_{r-1,1}\twoheadrightarrow\mathcal{Q}_{r-2}:={\rm coker}(\sK_{r-2}\hookrightarrow \mathcal{F}^{1}_{r-1,1}\oplus \mathcal{F}^{2}_{r-1,1})$, which induces $\varphi_1^{-1}$ on $U_0$.

When $1<j<\frac{r}{2}$, let $\sK_{r-2j}:={\rm ker}(\mathcal{F}^{1}_{r-j,j}\oplus \mathcal{F}^{2}_{r-j,j}\xrightarrow{q}\mathcal{Q}_{r-2j}\to 0)$ and
$$\mathcal{K}_{r-2(j-1)}:={\rm ker}(\mathcal{F}^{1,\vee}_{r-j,j-1}\oplus \mathcal{F}^{2,\vee}_{r-j,j-1}\twoheadrightarrow\mathcal{F}^{1,\vee}_{r-j,j}\oplus \mathcal{F}^{2,\vee}_{r-j,j}\twoheadrightarrow\mathcal{K}^{\vee}_{r-2j})$$
(use $0\to \mathcal{F}^i_{r-j,j}\to\mathcal{F}^i_{r-j,j-1}\to \mathcal{F}^i_{j,j-1}\to 0$).
Then, by the exact sequences $0\to\mathcal{F}^i_{r-(j-1),r-j}\to\mathcal{F}^i_{r-(j-1),j-1}\to \mathcal{F}^i_{r-j,j-1}\to 0$
and $$0\to \sK_{r-2(j-1)}\to \mathcal{F}^{1,\vee}_{r-j,j-1}\oplus \mathcal{F}^{2,\vee}_{r-j,j-1}\to \mathcal{K}^{\vee}_{r-2j}\to 0,$$
we have
$\mathcal{F}^{1}_{r-(j-1),j-1}\oplus \mathcal{F}^{2}_{r-(j-1),j-1}\twoheadrightarrow\mathcal{F}^{1}_{r-j,j-1}\oplus \mathcal{F}^{2}_{r-j,j-1}\twoheadrightarrow\mathcal{K}^{\vee}_{r-2(j-1)}$,
which induces the morphism $\varphi_j:Gr_j\to Gr_{j-1}$.

It is clear that $\varphi_j(Gr_j)=D^{[j-1]}_1\cap D^{[j-1]}_2$. To see that $$Gr_j\xrightarrow{\varphi_j} D^{[j-1]}_1\cap D^{[j-1]}_2$$ is birational, let
$E_i\twoheadrightarrow\cdots\twoheadrightarrow \sQ_{i,r-(j-1)}\twoheadrightarrow \sQ_{i,j-1}\twoheadrightarrow\cdots$ be
the flags on $D^{[j-1]}_1\cap D^{[j-1]}_2$ and $\mathcal{F}^{1}_{r-(j-1),j-1}\oplus \mathcal{F}^{2}_{r-(j-1),j-1}\xrightarrow{q}\mathcal{Q}_{r-2(j-1)}$ be the universal quotient. Let $U_{j-1}\subset D^{[j-1]}_1\cap D^{[j-1]}_2$ be the open set where $\mathcal{F}^{i}_{r-(j-1),j-1}\xrightarrow{q_i}\mathcal{Q}_{r-2(j-1)}$ has rank $r-2j+1$, we will construct flags
$$E_i\twoheadrightarrow\cdots\twoheadrightarrow \sQ_{i,r-(j-1)}\twoheadrightarrow\sQ_{i,r-j}\twoheadrightarrow\sQ_j\twoheadrightarrow\sQ_{i,j-1}\twoheadrightarrow\cdots$$
on $U_{j-1}$ (extended the given flags), which together with a quotient $$\sF^1_{r-j,j}\oplus\sF^2_{r-j,j}\twoheadrightarrow \sQ_{r-2j}$$
define $\varphi^{-1}_j$ on $U_{j-1}$. Let $\sF^i_{r-(j-1),r-j}:={\rm ker}(\mathcal{F}^{i}_{r-(j-1),j-1}\xrightarrow{q_i}\mathcal{Q}_{r-2(j-1)})$,
$${\rm Im}(q_i)\cong\frac{\mathcal{F}^{i}_{r-(j-1),j-1}}{\sF^i_{r-(j-1),r-j}}:=\mathcal{F}^i_{r-j,j-1}\subset
\sQ_{i,r-j}:=\frac{\sQ_{i,r-(j-1)}}{\sF^i_{r-(j-1),r-j}}.$$
Then $\mathcal{F}^{1}_{r-(j-1),j-1}\oplus \mathcal{F}^{2}_{r-(j-1),j-1}\xrightarrow{q}\mathcal{Q}_{r-2(j-1)}$ factors through
$$q: \mathcal{F}^{1}_{r-(j-1),j-1}\oplus \mathcal{F}^{2}_{r-(j-1),j-1}\twoheadrightarrow \mathcal{F}^{1}_{r-j,j-1}\oplus \mathcal{F}^{2}_{r-j,j-1}\twoheadrightarrow\mathcal{Q}_{r-2(j-1)}.$$
Let $\sK_{r-2j}:={\rm ker}(\mathcal{F}^{1}_{r-j,j-1}\oplus \mathcal{F}^{2}_{r-j,j-1}\twoheadrightarrow\mathcal{Q}_{r-2(j-1)})$, then homomorphisms
$\sK_{r-2j}\xrightarrow{p_i}\mathcal{F}^{i}_{r-j,j-1}$ (induced by projections $\mathcal{F}^{1}_{r-j,j-1}\oplus \mathcal{F}^{2}_{r-j,j-1}\to \mathcal{F}^{i}_{r-j,j-1}$) are
injective. Let $\sF^i_{r-j,j}:=p_i(\sK_{r-2j})\subset \mathcal{F}^{i}_{r-j,j-1}$,
$$\mathcal{Q}^{\vee}_{r-2j}:={\rm ker}(\mathcal{F}^{1,\vee}_{r-j,j}\oplus \mathcal{F}^{2,\vee}_{r-j,j}\twoheadrightarrow\mathcal{K}^{\vee}_{r-2j}),$$
$\sQ_{i,j}:=\frac{\sQ_{i,r-j}}{\sF^i_{r-j,j}}$. Then we have $\mathcal{F}^{1}_{r-(j-1),j}\oplus \mathcal{F}^{2}_{r-(j-1),j}\xrightarrow{q}\mathcal{Q}_{r-2j}\to 0$ and
$\mathcal{Q}_{i,r-(j-1)}\twoheadrightarrow\mathcal{Q}_{i,r-j}\twoheadrightarrow\mathcal{Q}_{i,j}\twoheadrightarrow\mathcal{Q}_{i,j-1}$, which
induce $\varphi^{-1}_j$ on $U_{j-1}$.

When $r$ is an even number and $j=\frac{r}{2}$, let
$$E_i\twoheadrightarrow\sQ_{i,r-1}\twoheadrightarrow\cdots\twoheadrightarrow\sQ_{i,\frac{r}{2}+1}\twoheadrightarrow\sQ_{i,\frac{r}{2}}\twoheadrightarrow\sQ_{i,\frac{r}{2}-1}\twoheadrightarrow\cdots$$
be universal flags on $Gr_{\frac{r}{2}}:=\mathrm{Flag}_{(1,\ldots, 1)}(E_1)\times_M \mathrm{Flag}_{(1,\ldots, 1)}(E_2)$,
$$\sF^i_{\frac{r}{2}+1,\frac{r}{2}-1}:={\rm ker}(\sQ_{i,\frac{r}{2}+1}\twoheadrightarrow\sQ_{i,\frac{r}{2}-1}),$$
and $\sF^i_{\frac{r}{2},\frac{r}{2}-1}:={\rm ker}(\sQ_{i,\frac{r}{2}}\twoheadrightarrow\sQ_{i,\frac{r}{2}-1})$. Then we have $\sF^i_{\frac{r}{2}+1,\frac{r}{2}-1}\twoheadrightarrow \sF^i_{\frac{r}{2},\frac{r}{2}-1}$,
and $\sF^1_{\frac{r}{2}+1,\frac{r}{2}-1}\oplus \sF^2_{\frac{r}{2}+1,\frac{r}{2}-1}\twoheadrightarrow \sF^1_{\frac{r}{2},\frac{r}{2}-1}\oplus \sF^2_{\frac{r}{2},\frac{r}{2}-1}$ defines the morphism
$$\varphi_{\frac{r}{2}}:Gr_{\frac{r}{2}}\to Gr_{\frac{r}{2}-1}:={\rm Grass}_2(\sF^1_{\frac{r}{2}+1,\frac{r}{2}-1}\oplus \sF^2_{\frac{r}{2}+1,\frac{r}{2}-1}).$$
It is easy to see that ${\rm Im}(\varphi_{\frac{r}{2}})=D_1^{[\frac{r}{2}-1]}\cap D_2^{[\frac{r}{2}-1]}$ and
$$\varphi_{\frac{r}{2}}:Gr_{\frac{r}{2}}\to D_1^{[\frac{r}{2}-1]}\cap D_2^{[\frac{r}{2}-1]}$$
is an isomorphism. Indeed, on $D_1^{[\frac{r}{2}-1]}\cap D_2^{[\frac{r}{2}-1]}$, there are flags
$$E_i\twoheadrightarrow\sQ_{i,r-1}\twoheadrightarrow\cdots\twoheadrightarrow\sQ_{i,\frac{r}{2}+1}\twoheadrightarrow\sQ_{i,\frac{r}{2}-1}\twoheadrightarrow\cdots$$
and universal quotient $\sF^1_{\frac{r}{2}+1,\frac{r}{2}-1}\oplus \sF^2_{\frac{r}{2}+1,\frac{r}{2}-1}\xrightarrow{q}\mathcal{Q}_2\to 0$ such that $\sF^i_{\frac{r}{2}+1,\frac{r}{2}}:={\rm ker}(\sF^i_{\frac{r}{2}+1,\frac{r}{2}-1}\xrightarrow{q_i}\sQ_2)$ are rank $1$ subbundles. The flags $$\cdots\twoheadrightarrow\sQ_{i,\frac{r}{2}+1}\twoheadrightarrow\sQ_{i,\frac{r}{2}}:=\sQ_{i,\frac{r}{2}+1}/\sF^i_{\frac{r}{2}+1,\frac{r}{2}}
\twoheadrightarrow\sQ_{i,\frac{r}{2}-1}\twoheadrightarrow\cdots$$
define the $\varphi_{\frac{r}{2}}^{-1}$ on $D_1^{[\frac{r}{2}-1]}\cap D_2^{[\frac{r}{2}-1]}$.
\end{proof}

For any $\sigma_j\in {\rm H}^0(Gr_j,\omega^{-1}_{Gr_j})$, we use the same symbol to denote its extension
$\sigma_j\in {\rm H}^0(\omega^{-1}_{D_1^{[j-1]}\cap D_2^{[j-1]}})$ (see the proof of Lemma \ref{lem4.1}).

\begin{lem}\label{lem4.3} When $p>3r$, for any line bundle $L$ on $\mathrm{Flag}_{\overrightarrow m_{j}}$,
\ga{4.6}{{\rm H}^0(\omega^{-1}_{Gr_{j}}(-D_1^{[j]}-D_2^{[j]})\otimes L)\to {\rm H}^0(\omega^{-1}_{D_1^{[j]}\cap D_2^{[j]}}\otimes L)}
$(0\le j<\frac{r}{2})$ are surjective. For $\sigma_j\in {\rm H}^0(Gr_j,\omega^{-1}_{Gr_j})$ $(1\le j\le\frac{r}{2})$, let
$\tilde{\sigma}_j\in {\rm H}^0(\omega^{-1}_{Gr_{j-1}}(-D_1^{[j-1]}-D_2^{[j-1]})$ such that $\tilde{\sigma}_j|_{D_1^{[j-1]}\cap D_2^{[j-1]}}=\sigma_j$, and \ga{4.7}{\sigma_{j-1}:=\tilde{\sigma}_j\cdot s_1^{[j-1]}\cdot s_2^{[j-1]}\in {\rm H}^0(Gr_{j-1},\omega^{-1}_{Gr_{j-1}})}
where $s_i^{[j-1]}\in {\rm H}^0(\sO_{Gr_{j-1}}(D_i^{[j-1]}))$ defines $D^{[j-1]}_i$. Then $\sigma_{j-1}$ vanishes on $D_1^{[j-1]}+D_2^{[j-1]}$ and
$\sigma^{p-1}_{j-1}$ splits $Gr_{j-1}$ if $\sigma^{p-1}_j$ splits $Gr_j$ and $M$ is a smooth variety with ${\rm H}^0(\sO_M)=k$.
.
\end{lem}
\begin{proof} We prove only surjectivity of \eqref{4.6} for $j=1$ (other cases are similar).
In this case, it is equivalent to prove that
\ga{4.8}{{\rm H}^0(Gr_0,\omega^{-1}_{Gr_0}(-D_1-D_2))\otimes L)\to {\rm H}^0(D_1\cap D_2,\omega^{-1}_{D_1\cap D_2}\otimes L)}
is surjective when $p>3r$. Let $\rho^*E_1\oplus \rho^*E_2\xrightarrow q\sQ\to 0$
be the universal quotient on $Gr_0\xrightarrow{\rho} M$. Then $\sO(D_i)=\rho^*{\rm det}(E_i)^{-1}\otimes{\rm det}(\sQ)$ and
$$\omega^{-1}_{Gr_0}(-D_1-D_2)={\rm det}(\sQ)^{2r-2}\otimes\rho^*(\omega^{-1}_M\otimes{\rm det}(E_1)^{1-r}\otimes {\rm det}(E_2)^{1-r}).$$
To prove that \eqref{4.8} is surjective, it is enough to show the surjection
\ga{4.9}{\rho_*({\rm det}(\sQ)^{2r-2})\to \rho^{+}_*({\rm det}(\sQ)^{2r-2}|_{D_1\cap D_2})}
has a splitting where $\rho^{+}:=\rho|_{D_1\cap D_2}:D_1\cap D_2\to M$. To see it, let $V_1$, $V_2$ be $r$-dimensional
vector spaces, $Gr:=\mathrm{Grass}_r(V_1\oplus V_2)$, $$V_1\otimes\sO_{Gr}\oplus V_2\otimes\sO_{Gr}\to \sQ\to 0$$
be the universal quotient bundle on $Gr$, ${\rm det}(\sQ):=\sO_{Gr}(1)$ and
$$S_i=\{\,V_1\oplus V_2\xrightarrow{q} Q\to 0\,|\, V_i\cap \mathrm{ker}(q)\neq 0\}\subset Gr,$$
then the restriction map
\ga{4.10} { {\rm H}^0(Gr, \sO_{Gr}(2r-2))\to {\rm H}^0(S_1\cap S_2, \sO_{S_1\cap S_2}(2r-2))} is a surjective homomorphism of ${\rm GL}(V_1)\times {\rm GL}(V_2)$-modules. In fact,
direct images in \eqref{4.9} can be thought of as the vector bundles associated to representations of  ${\rm GL}(V_1)\times {\rm GL}(V_2)$ in \eqref{4.10}, and we only need to show the homomorphism \eqref{4.10} of ${\rm GL}(V_1)\times {\rm GL}(V_2)$-modules
has a splitting, which holds in characteristic zero automatically.
Fortunately, when $p>3r$ it is also valid by the result from the representation of $\mathrm{GL}(r)$ in positive characteristic (see Lemma \ref{lem2.16}). We make other conclusions by following lemma.
\end{proof}

\begin{lem}\label{lem4.4} Let $\sM$ be a smooth variety with ${\rm H}^0(\sM,\sO_{\sM})=k$ and $D_i\subset \sM$ $(1\le i\le m)$ be prime divisors such that
$D_{j_1}\cap \cdots\cap D_{j_s}$ $(1\le s\le m)$ are normal varieties for any $1\le j_1<\cdots<j_s\le m$. Assume
$D_{j_1}\cap \cdots\cap D_{j_s}\subsetneqq D_{j_1}\cap\cdots\cap D_{j_{t-1}}\cap D_{j_{t+1}}\cap\cdots\cap D_{j_s}$ $(1\le t\le s)$.
Let $D_i$ be defined by $s_i\in {\rm H}^0(\sM, \sO_{\sM}(D_i))$ and $Z=D_1\cap \cdots\cap D_m$.
If $$\sigma\in {\rm H}^0(\sM, \omega^{-1}_{\sM}(-D_1-\cdots-D_m))$$ is a section such that $(\sigma|_{Z})^{p-1}$ splits $Z$, let
$$\wt\sigma:=\sigma\cdot s_1\cdots s_m\in {\rm H}^0(\sM, \omega^{-1}_{\sM}).$$
Then $\wt\sigma$ vanishes on $D_1+\cdots+D_m$ and ${\wt\sigma}^{p-1}$ splits $\sM$.
\end{lem}
\begin{proof} We need the following criteria (see \cite[Theorem 1.3.8]{M.S}): Let $\sM$ be a smooth variety with ${\rm H}^0(\sM,\sO_{\sM})=k$, $x_1,\cdots, x_n$ be a system of local coordinates at $y\in \sM$. Then $\varphi\in {\rm H}^0(\sM,\omega_{\sM}^{1-p})$ splits $\sM$ if and only if the monomial $(x_1\cdot \cdots \cdot x_n)^{p-1}$ occurs with coefficient $1$ in the local expansion of $\varphi$ at some (and hence every) point $y\in \sM$.

We prove the lemma by induction on $m$. When $m=1$, $$\sigma\in {\rm H}^0(\sM, \omega^{-1}_{\sM}(-D_1))$$ and $(\sigma|_{D_1})^{p-1}$ splits $D_1$, let $y\in D_1$ such that $\sM$ and $D_1$ are smooth at $y$, there exist local coordinates $x_1,x_2,...,x_n$ of $\sM$ at $y\in \sM$ such that $x_1=0$ is the local equation of $D_1$ at $y\in \sM$. Since $(\sigma|_{D_1})^{p-1}$ splits $D_1$, by the criteria, local expansion of $\sigma^{p-1}$ at $y$ contains monomial $(x_2\cdot \cdots \cdot x_n)^{p-1}$. Then the local expansion of $(\sigma\cdot s_1)^{p-1}$ at $y$ contains monomial $(x_1\cdot x_2 \cdots x_n)^{p-1}$. Thus $(\sigma\cdot s_1)^{p-1}$ splits $\sM$.

Assume the lemma holds for $m-1$. Let $D'_i:=D_1\cap D_i$, $\sigma':=\sigma|_{D_1}$ and $s'_i:=s_i|_{D_1}$ ($2\le i\le m$). Then $D_i'$ is a prime divisor of $D_1$ defined by $s_i'$ ($2\le i\le m$) and $\sigma'\in {\rm H}^0(D_1, \omega^{-1}_{D_1}(-D_2'-\cdots-D_m'))$ such that
$$\sigma'|_{D_2'\cap\cdots\cap D'_m}=\sigma|_Z,\quad Z=D_2'\cap\cdots\cap D'_m.$$
By the inductional assumption, let $\wt{\sigma}':=\sigma'\cdot s_2'\cdots s_m'\in {\rm H}^0(D_1, \omega^{-1}_{D_1})$, then $\wt{\sigma}^{'p-1}$ splits
$D_1$. Clearly, $\sigma\cdot s_2\cdots s_m\in {\rm H}^0(\sM, \omega^{-1}_{\sM}(-D_1))$ satisfies
$$(\sigma\cdot s_2\cdots s_m)|_{D_1}=\sigma'\cdot s_2'\cdots s_m'=\wt{\sigma}'\in {\rm H}^0(D_1, \omega^{-1}_{D_1}).$$
Thus, by using case $m=1$ of the lemma, we are done.
\end{proof}

To sum up, we have $Gr_{\frac{r}{2}}=\mathrm{Flag}_{(1,\ldots,1)}(E_1)\times_M \mathrm{Flag}_{(1,\ldots,1)}(E_2):=Y$ when $r$ is an even number, and when $r$ is an odd number,
$$Gr_{\frac{r-1}{2}}:={\rm Grass}_1(\sF^1_{\frac{r+1}{2},\frac{r-1}{2}}\oplus\sF^1_{\frac{r+1}{2},\frac{r-1}{2}})\xrightarrow{\rho_{\frac{r-1}{2}}}\mathrm{Flag}_{\overrightarrow m_{\frac{r-1}{2}}}=Y$$
is a $\mathbb{P}^1$-bundle over $Y$, $D^{[\frac{r-1}{2}]}_1$, $D^{[\frac{r-1}{2}]}_2$ are two disjoint sections. For any line bundle $L$ on $M$, when $p>3r$, the homomorphism
$${\rm H}^0(\omega^{-1}_{Gr_{\frac{r-1}{2}}}(-D^{[\frac{r-1}{2}]}_i)\otimes L)\to {\rm H}^0(\omega^{-1}_{D^{[\frac{r-1}{2}]}_i}\otimes L)$$
is surjective. For any $\sigma\in {\rm H}^0(\omega_Y^{-1}\otimes L)={\rm H}^0(\omega^{-1}_{D^{[\frac{r-1}{2}]}_1}\otimes L)$, we choose
\ga{4.11}{\tilde{\sigma}\in {\rm H}^0(\omega^{-1}_{Gr_{\frac{r-1}{2}}}(-D^{[\frac{r-1}{2}]}_1)\otimes L),\quad \tilde{\sigma}|_{D^{[\frac{r-1}{2}]}_1}=\sigma. }
Then we summarize our constructions in the following corollary.

\begin{cor}\label{cor4.5}
Let $E_1$, $E_2$ be vector bundles of rank $r$ on a variety $M$ and $Y:=\mathrm{Flag}_{(1,\ldots,1)}(E_1)\times_M \mathrm{Flag}_{(1,\ldots,1)}(E_2)$. When $p>3r$, for any $\sigma\in {\rm H}^0(\omega_Y^{-1}\otimes L)$, we construct a section $\sigma_0\in {\rm H}^0(\omega_{Gr_0}^{-1}\otimes L)$, where $Gr_0={\rm Grass}_r(E_1\oplus E_2)$, in the following way.
\begin{itemize}
\item [(1)] When $r$ is even, let $\sigma_{\frac{r}{2}}:=\sigma$ and, for any $1\le j\le\frac{r}{2}$,
$$\sigma_{j-1}:=\tilde{\sigma}_j\cdot s_1^{[j-1]}\cdot s_2^{[j-1]}\in {\rm H}^0(Gr_{j-1},\omega^{-1}_{Gr_{j-1}}\otimes L).$$
\item [(2)] When $r$ is odd, let $\sigma_{\frac{r-1}{2}}:=\tilde{\sigma}\cdot s_1^{[\frac{r-1}{2}]}\in {\rm H}^0(\omega^{-1}_{Gr_{\frac{r-1}{2}}}\otimes L)$ and
$$\sigma_{j-1}:=\tilde{\sigma}_j\cdot s_1^{[j-1]}\cdot s_2^{[j-1]}\in {\rm H}^0(Gr_{j-1},\omega^{-1}_{Gr_{j-1}}\otimes L)$$
for any $1\le j\le\frac{r}{2}-1$.
\item [(3)] If $\sigma_j^{p-1}$ splits $Gr_j$ and $M$ is smooth with ${\rm H}^0(\sO_M)=k$, then $\sigma_{j-1}$ vanishes on $D_1^{[j-1]}+D_2^{[j-1]}$ and $\sigma_{j-1}^{p-1}$ splits $Gr_{j-1}$. In particular, $\sigma_0$ vanishes on $D_1+D_2$ and $\sigma_0^{p-1}$ splits $Gr_0=\mathrm{Grass}_r(E_1\oplus E_2)$ if $\sigma^{p-1}$ splits
    $Y=\mathrm{Flag}_{(1,\ldots,1)}(E_1)\times_M \mathrm{Flag}_{(1,\ldots,1)}(E_2)$.
\end{itemize}
\end{cor}

\subsection{Proof of Theorem \ref{thm3.8}} Let $\omega_c=(2r,\{\vec {a}(x),\vec{n}(x)\}_{x\in I})$,
$$\wt{\sR}^{ss}_{\omega_c}\subset \wt{R}_F,\qquad \wt{\sR}^{ss}_{\omega_c}\to \wt{\sR}^{ss}_{\omega_c}//{\rm SL(V)}:=\sU_{\wt X,\,\omega_c}.$$
Let $\wt{\omega}_c=(2r,\{\vec{a}(x),\vec{n}(x)\}_{x\in I\cup\{x_1,\,x_2\}})$ be the data extending $\omega_c$ by adding points $x_1$, $x_2$ with type
$\vec{n}(x_i)=(1,\cdots,1)$,
$$Y:=\mathrm{Flag}_{\vec{n}(x_1)}(\sE_{x_1})\times_{\wt{R}_F} \mathrm{Flag}_{\vec{n}(x_2)}(\sE_{x_2}),$$
$Y^{ss}_{\wt\omega_c}\subset Y$ and $\psi:Y^{ss}_{\wt\omega_c}\to Y^{ss}_{\wt\omega_c}//{\rm SL}(V):=\sU_{\wt X,\,\wt\omega_c}$.

Let $\sigma\in {\rm H}^0(\omega^{-1}_{\sU_{\wt X,\,\wt\omega_c}})$ be a section such that $\sigma^{p-1}$ splits $\sU_{\wt X,\,\wt\omega_c}$ and
$$\hat{\sigma}:=\psi^*(\sigma)\in {\rm H}^0(Y^{ss}_{\wt\omega_c},\omega^{-1}_Y)^{inv.}= {\rm H}^0(Y,\omega^{-1}_Y)^{inv.}.$$

\begin{lem}\label{lem4.6} $\forall\,\,\hat{\sigma}\in {\rm H}^0(Y,\omega^{-1}_Y)^{inv.}$, there is a $\hat{\sigma}_0\in {\rm H}^0(\wt{\sR}'_F,\omega^{-1}_{\wt{\sR}'_F})^{inv.}$
which is obtained by constructions in Lemma \ref{lem4.3} and Corollary \ref{cor4.5}, where
$\wt{\sR}'_F={\rm Grass}_r(\sE_{x_1}\oplus\sE_{x_2})\xrightarrow{\rho}\wt\sR_F$.
\end{lem}

Lemma \ref{lem4.6} is always true in characteristic zero since surjectivity of ${\rm H}^0(\omega^{-1}_{Gr_{j}}(-D_1^{[j]}-D_2^{[j]})\otimes L)\twoheadrightarrow {\rm H}^0(\omega^{-1}_{D_1^{[j]}\cap D_2^{[j]}}\otimes L)$ implies the surjectivity of
${\rm H}^0(\omega^{-1}_{Gr_{j}}(-D_1^{[j]}-D_2^{[j]})\otimes L)^{inv.}\twoheadrightarrow {\rm H}^0(\omega^{-1}_{D_1^{[j]}\cap D_2^{[j]}}\otimes L)^{inv.}$
in this case. We do not know if it holds in general for positive characteristic. Recall $\wt{\sR}^{'ss}_{\omega_c,\alpha}\subset \wt{\sR}'_F$ and
$\wt{\sR}^{'ss}_{\omega_c,\alpha}\xrightarrow{\psi_{\alpha}} \sP_{\omega_c,\alpha}:=\wt{\sR}^{'ss}_{\omega_c,\alpha}//{\rm SL}(V)$, we have

\begin{lem}\label{lem4.7} Assume Lemma \ref{lem4.6} holds, let $\hat{\sigma}_0\in {\rm H}^0(\wt{\sR}'_F,\omega^{-1}_{\wt{\sR}'_F})^{inv.}$ be obtained from $\hat{\sigma}$
by constructions in Lemma \ref{lem4.3} and Corollary \ref{cor4.5}, and
\ga{4.12}{\sigma_{0,1^-}\in {\rm H}^0(\sP_{\omega_c,1^-}, \omega^{-1}_{\sP_{\omega_c,1^-}}),\quad \psi^*_{1^-}(\sigma_{0,1^-})=\hat{\sigma}_0|_{\wt{\sR}^{'ss}_{\omega_c,1^-}}.}
Then, if $\sigma^{p-1}$ splits $\sU_{\wt X,\,\wt\omega_c}$,
the section $\sigma_{0,1^-}^{p-1}$ splits $\sP_{\omega_c,1^-}$.
\end{lem}

\begin{proof} This is a local problem and we only need to show it for $\alpha=0^+$. We can also replace $\omega_c$ by a minor modification $\omega$ such that
there is a universal bundle $\sE$ on $\wt X\times \sU_{\wt X,\omega}$ and $\wt{\sR}^s_{\omega}=\wt{\sR}^{ss}_{\omega}$. Then
$$\sP_{\omega,0^+}={\rm Grass}_r(\sE_{x_1}\oplus\sE_{x_2})\xrightarrow{\rho}\sU_{\wt X,\omega}.$$
We can further choose $\vec {a}(x_i)$ ($i=1,\,2$) such that $Y^s_{\wt\omega}=Y^{ss}_{\wt\omega}$ and
$$\sU_{\wt X,\wt\omega}=\mathrm{Flag}_{\vec{n}(x_1)}(\sE_{x_1})\times_{\sU_{\wt X,\omega}} \mathrm{Flag}_{\vec{n}(x_2)}(\sE_{x_2}).$$

On the other hand, since $Y^s_{\wt\omega}=Y^{ss}_{\wt\omega}$, $\omega^{-1}_{Y^{ss}_{\wt\omega}}$ descents to $\omega^{-1}_{\sU_{\wt X,\wt\omega}}$
and $\hat{\sigma}|_{Y^{ss}_{\wt\omega}}$ descents to a section $\sigma'\in {\rm H}^0(\omega^{-1}_{\sU_{\wt X,\wt\omega}})$. Because $\sigma^{p-1}$ splits
$\sU_{\wt X,\wt\omega_c}$ and $\sigma=\sigma'$ on a common open set of $\sU_{\wt X,\wt\omega}$ and $\sU_{\wt X,\wt\omega_c}$, by local criteria of Frobenius splitting,
$(\sigma')^{p-1}$ splits $\sU_{\wt X,\wt\omega}$.

Now $\hat{\sigma}_0|_{\wt{\sR}^{'ss}_{\omega,0^+}}$ descents to a section $\sigma_{0,0^+}\in {\rm H}^0(\sP_{\omega,0^+}, \omega^{-1}_{\sP_{\omega,0^+}})$, which is
obtained from $\sigma'$ by constructions in Lemma \ref{lem4.3} and Corollary \ref{cor4.5}, thus $\sigma_{0,0^+}^{p-1}$ splits $\sP_{\omega,0^+}$. But $\sigma_{0,1^-}=\sigma_{0,0^+}$ on a common open set of $\sP_{\omega_c,1^-}$ and $\sP_{\omega,0^+}$, which implies that $\sigma_{0,1^-}^{p-1}$ splits $\sP_{\omega_c,1^-}$.
We are done.
\end{proof}

To complete the proof of Theorem \ref{thm3.8}, it is enough to finish the proof of Lemma \ref{lem4.6}. We will need the following remark.
\begin{rmk}\label{rmk4.8} For $\omega=(k,\{\vec n(x), \vec a(x)\}_{x\in I})$, let $\omega^{+}=(k,\{\vec n(x), \vec a(x)\}_{x\in I^+})$ where $I^+=I\cup\{z\}$, $\vec n(z)=(r-1,1)$ and $\vec a(z)=(0,1)$. Then
$$k\cdot{\rm par}_{\omega^+}\chi(E):=k\cdot\chi(E)+\sum_{x\in
I^+}\sum^{l_x+1}_{i=1}a_i(x)n_i(x)=k\cdot{\rm par}_{\omega}\chi(E)+1.$$
Thus, for any data $\omega=(k,\{\vec n(x), \vec a(x)\}_{x\in I})$, we can get a data $$\omega^+=(k,\{\vec n(x), \vec a(x)\}_{x\in I^+})$$ (by adding parabolic points)
such that $\omega^+$ satisfies
$$(r, k\cdot{\rm par}_{\omega^+}\chi(E))=1,\quad\underset{0<r_1<r}{\mathrm{min}}\{ r_1(r-r_1)(g-2)+\frac{|I^+|}{{\rm max}\{k,\,2r\}}\}\ge 2.$$
\end{rmk}

\begin{proof}[Proof of Lemma \ref{lem4.6}]
Let $\omega^+=(2r,\{\vec {a}(x),\vec{n}(x)\}_{x\in I^+})$ be the data extending $\omega_c=(2r,\{\vec {a}(x),\vec{n}(x)\}_{x\in I})$ by adding parabolic points, which
satisfies conditions in Remark \ref{rmk4.8}. Let $\sE^+:=({\rm id}\times\pi)^*\sE$ where $\sE$ is the universal quotient on $\wt X\times \wt\sR_F$,
$${\wt\sR}^+_F=\underset{x\in I^+}{\times_{\wt{\bold Q}_F}}{\rm Flag}_{\vec n(x)}(\wt\sF_x)\xrightarrow{\pi} \wt\sR_F=\underset{x\in I}{\times_{\wt{\bold Q}_F}}{\rm Flag}_{\vec n(x)}(\wt\sF_x) $$
and
$Y^+:=\mathrm{Flag}_{\vec{n}(x_1)}(\sE^+_{x_1})\times_{{\wt{R}}^+_F} \mathrm{Flag}_{\vec{n}(x_2)}(\sE^+_{x_2})\xrightarrow{\pi} Y$ are projections. Then
$\pi^*\omega^{-1}_Y=\omega^{-1}_{Y^+}\otimes L$ for some line bundle $L$ on ${\wt\sR}^+_F$. Let
$${\hat{\sigma}}^+:=\pi^*(\hat{\sigma})\in {\rm H}^0(\omega^{-1}_{Y^+}\otimes L)^{inv.}.$$
We \textbf{claim} that there is a ${\hat{\sigma}}^+_0\in {\rm H}^0(\omega^{-1}_{{\wt\sR}^{'+}_F}\otimes L)^{inv.}$, which is obtained from ${\hat{\sigma}}^+$ by constructions
in Lemma \ref{lem4.3} and Corollary \ref{cor4.5}, where
$${\wt\sR}^{'+}_F={\rm Grass}_r(\sE^+_{x_1}\oplus\sE^+_{x_2})\xrightarrow{\rho_+}{\wt\sR}^+_F, \quad \sE^+_{x_i}=\pi^*\sE_{x_i}$$
is obtained from ${\wt\sR}^{'}_F={\rm Grass}_r(\sE_{x_1}\oplus\sE_{x_2})\xrightarrow{\rho}{\wt\sR}_F$ by base change ${{\wt\sR}^+_F\xrightarrow{\pi}{\wt\sR}_F}$.
Assume the \textbf{claim}, let ${\wt\sR}^{'+}_F\xrightarrow{\pi_+}{\wt\sR}^{'}_F$ be the projection, then there is a $\hat{\sigma}_0\in {\rm H}^0(\omega^{-1}_{{\wt\sR}^{'}_F})^{inv.}$ such that ${\hat{\sigma}}^+_0=\pi_+^*(\hat{\sigma}_0)$ since
$${\rm H}^0(\omega^{-1}_{{\wt\sR}^{'}_F})\xrightarrow{\pi_+^*}{\rm H}^0(\pi_+^*\omega^{-1}_{{\wt\sR}^{'}_F})= {\rm H}^0(\omega^{-1}_{{\wt\sR}^{'+}_F}\otimes L)$$
is an isomorphism. It is easy to see that $\hat{\sigma}_0$ is obtained from $$\hat{\sigma}\in {\rm H}^0(\omega^{-1}_{Y})^{inv.}$$ by constructions
in Lemma \ref{lem4.3} and Corollary \ref{cor4.5} (note projections $Gr_j^+\xrightarrow{\pi_j} Gr_j$ induce isomorphisms ${\rm H}^0(\omega^{-1}_{Gr_j})\xrightarrow{\pi^*_j}
{\rm H}^0(\pi^*_j\omega^{-1}_{Gr_j})$).

Then it is enough to prove the \textbf{claim}. By choice of $\omega^+$, we have
$${\wt\sR}^{+,ss}_{\omega^+}={\wt\sR}^{+,s}_{\omega^+}\subset {\wt\sR}^+_F,\quad {\wt\sR}^{+,s}_{\omega^+}\xrightarrow{\psi_+}{\wt\sR}^{+,s}_{\omega^+}//\mathrm{SL}(V):=\sU_{\wt X,\,\omega^+}$$
and a universal object $(\mathcal{F},\, \{\mathcal{F}_x \twoheadrightarrow \mathcal{Q}_{l_x}\twoheadrightarrow\cdots \twoheadrightarrow \mathcal{Q}_{1}\twoheadrightarrow 0\}_{x\in I})$
on $\wt X\times\sU_{\wt X,\,\omega^+}$ such that $\sE^+=({\rm id}\times \psi_+)^*\sF$. Then ${\wt\sR}^{'+,ss}_{\omega^+,0^+}={\wt\sR}^{'+,s}_{\omega^+,0^+}=\rho_+^{-1}({\wt\sR}^{+,s}_{\omega^+})\subset{\wt\sR}^{'+}_F$ and the moduli space
$\sP_{\omega^+,\,0^+}:={\wt\sR}^{'+,s}_{\omega^+,0^+}//\mathrm{SL}(V)$ is
$$\sP_{\omega^+,\,0^+}={\rm Grass}_r(\sF_{x_1}\oplus\sF_{x_2})\to \sU_{\wt X,\,\omega^+}.$$

Let $\wt\omega=\omega^+\cup (k',\{\vec n(x), \vec a(x)\}_{x\in \{x_1,x_2\}})$ such that
$\vec n(x)=(1,1,...,1)$ and $\frac{a_i(x)}{k'}=0^+$ for $x\in \{x_1,x_2\}$. Then
$$Y^{+,ss}_{\wt\omega}=Y^{+,s}_{\wt\omega}=\mathrm{Flag}_{\vec{n}(x_1)}(\sE^+_{x_1})\times_{{\wt{R}}^{+,s}_{\omega^+}} \mathrm{Flag}_{\vec{n}(x_2)}(\sE^+_{x_2})\subset Y^+$$
and $Y^{+,s}_{\wt\omega}//\mathrm{SL}(V):=\mathcal{U}_{\widetilde{X},\,\wt\omega}=\mathrm{Flag}_{\vec{n}(x_1)}(\sF_{x_1})\times_{\sU_{\wt X,\,\omega^+}} \mathrm{Flag}_{\vec{n}(x_2)}(\sF_{x_2})$. Let $Y^{+,s}_{\wt\omega}\xrightarrow{\wt\psi_+}\mathcal{U}_{\widetilde{X},\,\wt\omega}$ be the quotient map.
Then it is clear that there is a line bundle $\sL$ on $\sU_{\wt X,\,\omega^+}$ and a section
$\sigma^+\in {\rm H}^0(\omega^{-1}_{\sU_{\wt X,\,\wt\omega}}\otimes \sL)$
such that $$(\omega^{-1}_{Y^+}\otimes L)|_{Y^{+,s}_{\wt\omega}}={\wt\psi_+}^*(\omega^{-1}_{\sU_{\wt X,\,\wt\omega}}\otimes \sL),\quad
{\hat{\sigma}}^+|_{Y^{+,s}_{\wt\omega}}={\wt\psi_+}^*(\sigma^+).$$
Apply constructions in Lemma \ref{lem4.3} and Corollary \ref{cor4.5} to $\sigma^+$, we get a section
$\sigma_0^+\in {\rm H}^0(\omega^{-1}_{\sP_{\omega^+,\,0^+}}\otimes\sL)$. Let ${\wt\sR}^{'+,s}_{\omega^+,0^+}\xrightarrow{\psi'_+}\sP_{\omega^+,\,0^+}$ be the quotient map.
Then ${\psi'_+}^*(\sigma_0^+)$ extends to a section ${\hat{\sigma}}^+_0\in {\rm H}^0(\omega^{-1}_{{\wt\sR}^{'+}_F}\otimes L)^{inv.}$ since
codimension of ${\wt\sR}^{'+}_F\setminus {\wt\sR}^{'+,s}_{\omega^+,0^+}$ has at least two. It is easy to see that ${\hat{\sigma}}^+_0$ is obtained from ${\hat{\sigma}}^+$ by constructions in Lemma \ref{lem4.3} and Corollary \ref{cor4.5} since the codimension of ${\wt\sR}^+_F\setminus{\wt\sR}^{+,s}_{\omega^+}$ is at least two. We are done.
\end{proof}

\subsection{Proof of Theorem \ref{thm3.17}} We only sketch it, which is similar with the proof of Theorem \ref{thm3.8}. Recall that
$$\sR_i=\underset{x\in I_i}{\times_{\textbf{Q}_i}}
{\rm Flag}_{\vec n(x)}(\sF_{i,x}),\quad \widetilde{\sR}={\rm Grass}_r(\sF_{x_1}\oplus\sF_{x_2})\xrightarrow{\rho}\sR:=\sR_1\times\sR_2$$
where $\sF:=\sF_1\oplus\sF_2$ (we use the same symbol $\sF_i$ to denote the universal quotient and its pullback on $X_i\times\sR_i$ and $\wt X\times\sR$), and
$$\wt\sR^{ss}_{\omega_c,\alpha}\xrightarrow{\psi_{\omega_c,\alpha}}\wt\sR^{ss}_{\omega_c,\alpha}//G:=\sP_{\omega_c,\alpha},\quad \sP_{\omega_c}:=\sP_{\omega_c,\,1}$$
where $\omega_c=(2r,\{\vec n(x),\vec a_c(x)\}_{x\in I_1\cup I_2})$. Let ${\rm Flag}_{\vec n(x_i)}(\sF_{x_i})^{ss}_{\omega_{c_i}}\xrightarrow{\psi_i}\sU_{X_i,\omega_{c_i}}$
be the quotient map, where ${\rm Flag}_{\vec n(x_i)}(\sF_{x_i})^{ss}_{\omega_{c_i}}\subset {\rm Flag}_{\vec n(x_i)}(\sF_{x_i})(\to \sR_i)$
is the open set of GIT semistable points respect to polarization $\omega^{-1}_{{\rm Flag}_{\vec n(x_i)}(\sF_{x_i})}$,
and $\theta_i\in {\rm H}^0(\omega^{-1}_{\sU_{X_i,\omega_{c_i}}})$ such that $\theta_i^{p-1}$ splits $\sU_{X_i,\omega_{c_i}}$. Then $\psi_i^*(\theta_i)$ can be extended to
a ${\rm GL}(V_i)$-invariant section $\hat{\theta}_i\in {\rm H}^0(\omega^{-1}_{{\rm Flag}_{\vec n(x_i)}(\sF_{x_i})})^{inv.}$. Let
$$\hat{\sigma}:=\hat{\theta}_1\boxtimes\hat{\theta}_2\in {\rm H}^0(\omega^{-1}_Y),\quad Y:={\rm Flag}_{\vec n(x_1)}(\sF_{x_1})\times_{\sR_F}{\rm Flag}_{\vec n(x_2)}(\sF_{x_2}).$$
Then $\hat{\sigma}$ is $G$-invariant where $G=(\mathrm{GL}(V_1)\times \mathrm{GL}(V_2))\cap \mathrm{SL}(V_1\oplus V_2)$. By a similar result of Lemma \ref{lem4.6},
we obtain a $G$-invariant section $\hat{\sigma}_0\in {\rm H}^0(\omega^{-1}_{\wt\sR_F})$ by constructions in Lemma \ref{lem4.3} and Corollary \ref{cor4.5}. Thus
$\hat{\sigma}_0|_{\wt\sR^{s}_{\omega_c,1^{-}}}$ descents to a section $\sigma_0\in {\rm H}^0(\omega^{-1}_{\sP_{\omega_c,1^-}})$. In a similar way of Lemma \ref{lem4.7}, we can show that $\sigma^{p-1}_0$ splits
$\sP_{\omega_c,1^-}$, which implies that there is a $\tilde{\sigma}\in {\rm H}^0(\omega^{-1}_{\sP_{\omega_c}})$ such that $\tilde{\sigma}$ vanishes on $D_1+D_2$ and $\tilde{\sigma}^{p-1}$ splits $\sP_{\omega_c}$ by Lemma \ref{lem4.1}. we are done.

\section{Moduli space of parabolic bundles on $\mathbb{P}^1$}
For any set $I$ of points with given type $\vec{n}(x)$ ($x\in I$), let $$\omega_c=(2r,\{\vec n(x),\,\vec a_c(x)\}_{x\in I})$$
be the canonical weight. The goal of this section is to prove
\begin{thm}\label{thm5.1}
Let $I_1=\{y_1,z_1,z_2,z\}$ be any set of points on $\mathbb{P}^1$, and
\ga{5.1} {\omega^1_c:=(2r,\{\vec n(x),\vec a_c(x)\}_{x\in I_1})}
with $\vec n(y_1)=(r-1,1)$, $\vec n(x)=(1,\ldots,1)$ for $x\in I-\{y_1\}$. Then $\mathcal{U}_{\mathbb{P}^1,\,\omega^1_c}$ is split by a $(p-1)$-power.
\end{thm}

We start with some general discussions. Our problem can be simplified to the case $d=0$ by Hecke transformation.
Given a parabolic sheaf $E$ with quasi-parabolic structure
$$E_z=Q_{r}(E)_z\twoheadrightarrow Q_{r-1}(E)_z\twoheadrightarrow \cdots \cdots \twoheadrightarrow Q_1(E)_z\twoheadrightarrow Q_0(E)_z=0$$ for a fixed point $z\in I$ of type $\vec n(z)=(1,\ldots, 1)$ with weight $$0=a_1(z)<a_2(z)<\ldots  <a_{r}(z)<k.$$ Let $F_i(E)_z=\mathrm{\mathrm{\mathrm{ker}}}\{E_z\twoheadrightarrow Q_i(E)_z\}$ and $E'=\mathrm{\mathrm{\mathrm{ker}}}\{E\twoheadrightarrow Q_1(E)_z\}$. Then at $z\in I$, $E'$ has a natural quasi-parabolic structure $$E_z'\twoheadrightarrow F_1(E)_z\twoheadrightarrow Q_{r-2}(E')_z \twoheadrightarrow \cdots \twoheadrightarrow Q_1(E')_z\twoheadrightarrow 0$$ of type ${\vec n}'(z)=(1,\ldots, 1)$, where $Q_i(E')_z\subset Q_{i+1}(E)_z$ ($1\leq i\leq r-2$) is the image of $F_1(E)_z$ under $E_z\twoheadrightarrow Q_{i+1}(E)_z$.

\begin{defn}\label{defn5.2} The parabolic sheaf $E'$ with given weight $$0=a'_1(z)<\cdots < a'_r(z)<k$$
is called Hecke transformation of parabolic sheaf $E$ at $z\in I$, where $a'_r(z)=k-a_2(z)$ and
$a'_i(z)=a_{i+1}(z)-a_2(z)$ for $2\leq i\leq r-1$ (Clearly, canonical weight goes to canoincal weght under Hecke transformation).
\end{defn}
One can also define the Hecke transformation of a family of parabolic sheaves (flat family yielding flat family, and preserve semistability). Thus, we have a morphism $$H_z: \sU_{C,\,\omega}=\sU_C(r,d,\omega)\rightarrow \sU_C(r,d-1,\omega')=\sU_{C,\,\omega'}.$$
\begin{lem}\label{lem5.3}
The morphism  $H_z$ is an isomorphism.
\end{lem}
\begin{proof} $H_z^{-1}$ was construced in the proof of \cite[Lemma 3.9]{SZ}.
\end{proof}

Let $\omega^0_c:=(2r,\{\vec n(x),\vec a_c(x)\}_{x\in \{y_1,z_1,z_2\}})$ be canonical weight and
\ga{5.2} {\omega=\omega^0_c\cup \omega':=(k,\{\vec n(x), \vec a(x)\}_{x\in\{y_1,z_1,z_2\}\cup I'})}
where $\omega'=(k,\{\vec n(x), \vec a(x)\}_{x\in I'})$ is an arbitrary weight and $\vec a(x)=(a_1(x)=0,a_2(x),\ldots,a_{l_x+1})$ satisfies $$\frac{a_{i+1}(x)-a_i(x)}{k}=\frac{n_i(x)+n_{i+1}(x)}{2r},\,\,\,\,\forall \,x\in \{y_1,z_1,z_2\}$$ then we have

\begin{prop}\label{prop5.4} $\mathcal{U}_{\mathbb{P}^1,\,\omega^0_c}$ consists of a $\omega^0_c$-stable parabolic bundle with underlying
bundle $E=\sO^{\oplus r}_{\mathbb{P}^1}$, and $\mathcal{U}_{\mathbb{P}^1,\,\omega}$ is birationally equivalent to
$$\underset{x\in\, I'}{\times}{\rm Flag}_{\vec n(x)}(E_x).$$
In particular, $\mathcal{U}_{\mathbb{P}^1,\,\omega^1_c}$ is birationally equivalent to ${X}:=\mathrm{Flag}_{(1,\ldots, 1)}(W)$
where $W=E_z$ is a vector space of dimension $r$.
\end{prop}
\begin{proof} Let $\sR_0:=\underset{x\in\{y_1,z_1,z_2\}}{\times_{\mathbf{Q}_F}}{\rm Flag}_{\vec n(x)}(\sF_x)$ and $\sR^{s}_{0,\omega_c^0}\subset \sR_{0}$ be the open set
of $\omega_c^0$-stable points. Then ${\rm Codim}(\sR_{0}\setminus\sR^{s}_{0,\omega_c^0})>0$ (see Remark \ref{rmk2.10}) and
$\mathcal{U}_{\mathbb{P}^1,\,\omega^0_c}:=\sR^{ss}_{0,\omega_c^0}//{\rm SL}(V)$ has dimension zero. Thus $\sR^{ss}_{0,\omega_c^0}=\sR^{s}_{0,\omega_c^0}$ and
$\mathcal{U}_{\mathbb{P}^1,\,\omega^0_c}$ consists of a $\omega^0_c$-stable parabolic bundle with underlying
bundle $E=\sO^{\oplus r}_{\mathbb{P}^1}$ (since $\sR_0$ is irreducible and the set of quotients with trivial underlying bundles is a nonempty open set of $\sR_0$). Let
$$\sR=\sR_0\times_{\mathbf{Q}_F}\left(\underset{x\in I'}{\times_{\mathbf{Q}_F}}{\rm Flag}_{\vec n(x)}(\sF_x)\right)\xrightarrow{\pi}\sR_0,\quad \sR^{ss}_{\omega}\xrightarrow{\psi}\sR^{ss}_{\omega}//{\rm SL}(V):=\sU_{\mathbb{P}^1,\,\omega},$$
$\hat{U}:=\pi^{-1}(\sR^s_{0,\omega^0_c})\cap \sR^s_{\omega}\subset \sR^s_{\omega}$, $U:=\psi(\hat{U})\subset \sU_{\mathbb{P}^1,\,\omega}$ and
$$\hat{U}_q:=\hat{U}\cap\pi^{-1}(q)\subset \underset{x\in\, I'}{\times}{\rm Flag}_{\vec n(x)}(E_x)$$
for $q=(V\otimes\sO_{\mathbb{P}^1}(-N)\twoheadrightarrow E,\, Q_{\bullet}(E)_{y_1},\, Q_{\bullet}(E)_{z_1},\,Q_{\bullet}(E)_{z_2})\in \sR^s_{0,\omega^0_c}$.

For any $q\in \sR^s_{0,\omega^0_c}$, the set $\hat{U}_q$ is not empty since $\sR^s_{0,\omega_c^0}$ is transitive under ${\rm SL}(V)$, and
\ga{5.3} {\psi_q:=\psi|_{\,\hat{U}_q}: \hat{U}_q\to U}
is an isomorphism. We show firstly that $\psi_q$ is bijective. For any two points $x=(q, u),\, y=(q,v) \in \hat{U}_q$, if
$\psi_q(x)=\psi_q(y)$, there is a $\eta\in {\rm SL}(V)$ such that $\eta\cdot x=y$. It means $\eta\cdot q=q$, which implies $\eta=\lambda\cdot{\rm id}_V$. Thus $y=\eta\cdot x=x$ and $\psi_q$ is injective. For any point $x=(q', v')\in \hat{U}=\pi^{-1}(\sR^s_{0,\omega^0_c})\cap \sR^s_{\omega}$, there is a $\eta\in {\rm SL}(V)$ such that $\eta\cdot q'=q$ since
$\sR^s_{0,\omega^0_c}//{\rm SL}(V)$ consists of one point.
Then $\eta\cdot x=(\eta\cdot q',\eta\cdot v')\in \hat{U}_q$ and $\psi(x)=\psi_q(\eta\cdot x)\in \psi_q(\hat{U}_q)$, thus $\psi_q$ is surjective. To
show that $\psi_q$ is an isomorphism, it is enough to check surjectivity of $\psi_q^*\Omega^1_{U}\to \Omega^1_{\hat{U}_q}$, which is true if
$\psi_q$ is \textit{formally unramified} by \cite[Corollary 17.2.2]{EGA1}. According to \cite[Definition 17.1.1]{EGA1}, $\psi_q:\hat{U}_q\to U$ is \textit{formally unramified} if, for any $X_0:={\rm Spec}(A/I)\subset X:={\rm Spec}(A)$ where $I^2=0$, the map
\ga{5.4}{{\rm Hom}_{U}(X, \hat{U}_q)\to {\rm Hom}_{U}(X_0, \hat{U}_q),\quad \varphi\mapsto \varphi|_{X_0}}
is injective. Since any $\varphi\in {\rm Hom}_{U}(X, \hat{U}_q)$ is determined by a set of flags $\{E_x\otimes A\twoheadrightarrow Q_{\bullet}(E)_x\}_{x\in I'}$ of finite $A$-modules, the injectivity of \eqref{5.4} follows from the Nakayama's lemma.
\end{proof}

\begin{rmk}\label{rmk5.5}
(1) When $I'=\{z\}$ and $r=2$, $\sU_{\mathbb{P}^1,\,\omega}=\mathbb{P}^1$ for any choice of $\frac{a_i(z)}{k}$. In particular, Theorem \ref{thm5.1} holds for $r=2$.

(2) If there is a point $q\in \sR^s_{0,\omega^0_c}$ (thus, for any point) such that
$\pi^{-1}(q)\subset\sR^{ss}_{\omega}$ and $\pi^{-1}(q)\cap\sR^{s}_{\omega}\neq \emptyset$, then \eqref{5.3} induces
$$\psi_q: \pi^{-1}(q)=\underset{x\in\, I'}{\times}{\rm Flag}_{\vec n(x)}(E_x)\cong \sU_{\mathbb{P}^1,\,\omega}.$$
Indeed, $\psi_q$ is birational by Proposition \ref{prop5.4}. But $\psi_q^*\Theta_{\sU_{\mathbb{P}^1,\,\omega}}=\Theta_{\sR,\,\omega}|_{\pi^{-1}(q)}$ is
ample, thus $\psi_q$ must be an isomorphism.
\end{rmk}

Assume $\vec n(x)=(1,\ldots,1)$ ($\forall\,\,x\in I'$) and $I=\{y_1,z_1,z_2\}\cup I'$. Let
\ga{5.5}{\omega^t=\omega_c^0\cup (k,\{\vec n(x),\,\vec a^t(x)\}_{x\in I'}),\quad |t|<\frac{1}{|I'|(r-1)r},}
where
$\frac{a^t_i(x)}{k}=\frac{i-1}{|I'|(r-1)r}$ ($\forall\,\,x\in I'\setminus\{z\}$, $1\le i\le r$), and
$$\frac{a^t_i(z)}{k}=\frac{i-1}{|I'|(r-1)r},\quad \frac{a^t_r(z)}{k}=\frac{1}{|I'|r}+t,\quad 1\le i\le r-1.$$
Let $t=0^+>0$ be sufficiently small such that $\sR^s_{\omega^{0^+}}=\sR^{ss}_{\omega^{0^+}}$. Then

\begin{lem}\label{lem5.6} (1) If $\sU_{\mathbb{P}^1,\,\omega^{0^+}}$ is split by a $(p-1)$-power, so is $\sU_{\mathbb{P}^1,\,\omega_c}$.

(2) When $t\le 0$, we have $\pi^{-1}(\sR^s_{0,\omega^0_c})\subset \sR^{ss}_{\omega^t}$. In particular, for any
$q\in \sR^s_{0,\omega^0_c}$,  $\sR^{ss}_{\omega^t}\xrightarrow{\psi}\sR^{ss}_{\omega^t}//{\rm SL}(V):=\sU_{\mathbb{P}^1,\,\omega^t}$ induces an isomorphism
\ga{5.6}{\psi_q: \pi^{-1}(q)=\underset{x\in\, I'}{\times}{\rm Flag}_{\vec n(x)}(E_x)\cong \sU_{\mathbb{P}^1,\,\omega^t}.}
\end{lem}

\begin{proof} (1) When $r=2$ and $|I'|=1$, $\sU_{\mathbb{P}^1,\,\omega^{0^+}}=\mathbb{P}^1=\sU_{\mathbb{P}^1,\,\omega_c}$ (see Remark \ref{rmk5.5}).
Thus we can assume either $r>2$ or $|I'|>1$. In this case, we have ${\rm Codim}(\sR\setminus \sR^s_{\omega_c})\ge 2$ (see Remark \ref{rmk2.10}). Thus it is enough to show
${\rm Codim}(\sR\setminus \sR^s_{\omega^{0^+}})\ge 2$. Indeed, note $\sR^s_{\omega^{0^+}}=\sR^{ss}_{\omega^{0^+}}$, we have
$$\mathrm{Codim}(\sR\setminus \sR^{s}_{\omega^{0^+}})>\underset{0<r_1<r}{\mathrm{min}}\left\{\mathrm{min}\{\frac{r-r_1}{2},\frac{r_1}{2}\} +\sum_{x\in I'}\Sigma_{x,r_1}(\omega^{0^+})\right\}$$(see Remark \ref{rmk2.10}). Thus, when $1<r_1<r-1$, we have
$$\mathrm{Codim}(\sR\setminus \sR^{ss}_{\omega^{0^+}})\ge 2.$$
When $r_1=1$, there is a $l$ such that $m_{l}(z)=1$, $m_i(z)=0$ ($i\neq l$),
$$\Sigma_{z,1}(\omega^{0^0})=l-1+\sum_{i\neq l}\frac{a_i^{0^+}(z)}{k}+(1-r)\frac{a_{l}^{0^+}(z)}{k}
\geq \sum_{1\leq i\leq r}\frac{a_i^{0^+}(z)}{k}\geq
\frac{1}{2}+0^+.$$
When $r_1=r-1$, there is a $l$ such that $m_{l}(z)=0$, $m_i(z)=1$ ($i\neq l$),
$$\Sigma_{z,r-1}(\omega^{0^+})=r-l-\sum_{i\neq l}\frac{a^{0^+}_{i}(z)}{k}+(r-1)\frac{a^{0^+}_{l}(z)}{k}
\geq \sum_{1\leq i\leq r}\frac{a_i^{0^+}(z)}{k}\geq
\frac{1}{2}+0^+.$$ Altogether we have $\mathrm{Codim}(\sR\setminus \sR^{s}_{\omega^{0^+}})\ge 2$.

(2) Let $Z=\sR\setminus\sR^{ss}_{\omega^t}$ and $U_0=\sR_0\setminus\pi(Z)$. If $U_0\neq \emptyset$, then
$$\pi^{-1}(U_0\cap \sR^s_{0,\omega_c^0})\subset \sR^{ss}_{\omega^t}$$
is not an empty  open subset. Thus by Remark \ref{rmk5.5} (2), we have the isomorphism \eqref{5.6} for any $q\in \sR^s_{0,\omega_c^0}$.

To show $U_0\neq \emptyset$, it is enough to show that ${\rm Codim}(\sR_0,\pi(Z))$
(the codimension of $\pi(Z)$ in $\sR_0$) is positive. Indeed
$$\aligned &\mathrm{Codim}(\sR_0,\pi(Z))>\\&\underset{0<r_1<r}{\rm min}\left\{\mathrm{min}\{\frac{r-r_1}{2},\frac{r_1}{2}\}+\sum_{x\in I'}\sum_{j=1}^{l_x+1}(r_1n_j(x)-rm_j(x))\frac{ a^t_j(x)}{k}\right\} \\&\geq
\underset{0<r_1<r}{\rm min}\left\{ \mathrm{min}\{\frac{r-r_1}{2},\frac{r_1}{2}\}+\frac{r_1}{2}+r_1t-rt-\frac{r_1(2r-r_1-1)}{2(r-1)}\right\}\\&\geq \underset{0<r_1<r}{\rm min}\left\{    (r_1-r)t      \right\}\ge 0 \,\,\, \text{when}\,\, t\le 0.
\endaligned$$
\end{proof}

In order to prove that $\sU_{\mathbb{P}^1,\,\omega^{0^+}}$ is split by a $(p-1)$-power when ${I'=\{\,z\,\}}$ (which implies Theorem \ref{thm5.1} by (1) of Lemma \ref{lem5.6}), we consider the birational morphism (induced by $\sR^s_{\omega^{0^+}}\subset \sR^{ss}_{\omega^{0}}$)
\ga{5.7}{f:\mathcal{U}_{\mathbb{P}^1,\,\omega^{0^+}}\rightarrow (\mathcal{U}_{
\mathbb{P}^1,\,\omega^0}\xrightarrow{\psi_q^{-1}})\mathrm{Flag}_{\vec n(z)}(E_z)}
(for any fixed $q\in \sR^s_{0,\,\omega_c^0}$). Since $f$ is an isomorphism on the open set
$$\psi(\pi^{-1}(\sR^s_{0,\omega^0_c})\cap \sR^{s}_{\omega^{0^+}})\subset \mathcal{U}_{\mathbb{P}^1,\,\omega^{0^+}}$$
(where $\sR^s_{\omega^{0^+}}\xrightarrow{\psi}\sR^s_{\omega^{0^+}}//{\rm SL}(V):=\sU_{\mathbb{P}^1,\,\omega^{0^+}}$), the exceptional divisors of $f$ must be located in
the following divisor
\ga{5.8} {\mathfrak{D}:=\psi(\pi^{-1}(\sR_0\setminus\sR^s_{0,\omega^0_c})\cap \sR^{s}_{\omega^{0^+}})\subset \mathcal{U}_{\mathbb{P}^1,\,\omega^{0^+}}.}
In fact, we will show that $\mathfrak{D}$ is the exceptional divisors of $f$ and consisting of $2r$ prime divisors. Recall $\sR^s_{\omega^{0^+}}\subset \sR^{ss}_{\omega^{0}}\subset \sR$ and
\ga{5.9}{\sR=\sR_0\times_{\mathbf{Q}_F}{\rm Flag}_{\vec n(z)}(\sF_z)\xrightarrow{\pi}\sR_0=\underset{x\in \{y_1,z_1,z_2\}}{\times_{\bold{Q}_F}}\mathrm{Flag}_{\vec n(x)}(\mathcal{F}_x).}
For any point $(V\otimes\sO_{\mathbb{P}^1}(-N)\twoheadrightarrow E,\, Q_{\bullet}(E)_{y_1},\, Q_{\bullet}(E)_{z_1},\,Q_{\bullet}(E)_{z_2},\,Q_{\bullet}(E)_z)$ of $\sR^{ss}_{\omega^{0}}$, we have either $E=\mathcal{O}_{\mathbb{P}^1}(-1)\oplus \mathcal{O}_{\mathbb{P}^1}(1)\oplus\mathcal{O}_{\mathbb{P}^1}^{\oplus (r-2)}$ or $E=\mathcal{O}_{\mathbb{P}^1}^{\oplus r}$ by
its parabolic semi-stability. Without loss of generality, we assume all points of $\sR$ have underlying bundle either ${\mathcal{O}_{\mathbb{P}^1}(-1)\oplus \mathcal{O}_{\mathbb{P}^1}(1)\oplus\mathcal{O}_{\mathbb{P}^1}^{\oplus (r-2)}}$ or $\mathcal{O}_{\mathbb{P}^1}^{\oplus r}$. Then
we will show
that $\sR_0\setminus\sR^s_{0,\omega^0_c}$ consists of exactly $2r$ prime divisors (using the fact that $\sR^s_{0,\omega^0_c}$ is an orbit).

To define the $2r$ prime divisors, let $V\otimes \sO_{\mathbb{P}^1\times \sR_0}(-N)\rightarrow \sE\rightarrow 0$ be the universal quotient, $0\subset\sE_{y_1,1}\subset \sE_{y_1}:=\sE|_{\{y_1\}\times \sR_0}$ and
$$0\subset \sE_{z_i,r-1}\subset\cdots\subset\sE_{z_i,2}\subset\sE_{z_i,1}\subset\sE_{z_i}:=\sE|_{\{z_i\}\times \sR_0}$$
be the universal flags. Then $pr_{2*}\sE$ is locally free of rank $r$ and commutes with base change (where
$pr_2:\mathbb{P}^1\times \sR_0\to \sR_0$ is the second projection),
$$\sR'_0:=\{\,q\in \sR_0\,|\,\sE|_{\mathbb{P}^1\times\{q\}}\cong \sO_{\mathbb{P}^1}^{\oplus r}\,\}\subset \sR_0$$ is the open subset
where the canonical homomorphism $pr_2^*pr_{2*}\sE\to \sE$ is an isomorphism. Thus $\sR_0\setminus \sR'_0$ is a prime
divisor of $\sR_0$. Let $\sE':=\sE|_{\mathbb{P}^1\times \sR'_0}$. The evaluation ${\rm H}^0(E)\xrightarrow{ev_x} E_x$ induces isomorphisms
$pr_{2*}\sE' \xrightarrow{ev_x}\sE'_x,\,\, \sE'_x\xrightarrow{ev_{xy}}\sE'_y,\,\,\,\, ev_{xy}:=ev_y\circ ev_x^{-1}$ for $x,y\in \{y_1,z_1,z_2\}$.
Then we firstly define the $2r$ prime divisors $\widetilde{\mathcal{X}}_i,\,\,\widetilde{\mathcal{Y}}_i\subset\sR_0$ ($1\leq i\leq r$).

Let $\widetilde{\mathcal{X}}_i\subset \mathcal{R}_0$ ($1\leq i\leq r$) be the closed subvariety whose generic point is determined by the condition that \ga{5.10}{\sE'_{y_1,1}\hookrightarrow \sE'_{y_1}\xrightarrow{ev^{-1}_{y_1}}pr_{2*}\sE'/(ev_{z_1}^{-1}(\sE'_{z_1,r-i+1})\oplus ev_{z_2}^{-1}(\sE'_{z_2,i}))} is a zero morphism, $\widetilde{\mathcal{Y}}_i \subset \sR_0$ ($1\leq i\leq r-1$) be the closed subvariety whose generic point is determined by the condition that \ga{5.11}{\sE'_{z_1,r-i}\hookrightarrow \sE'_{z_1}\xrightarrow{ev_{z_1}^{-1}} pr_{2*}\sE'/ev^{-1}_{z_2}(\sE'_{z_2,i})} is not an isomorphism, and $\widetilde{\mathcal{Y}}_r:=\sR_0\setminus \sR'_0$. Then we have

\begin{lem}\label{lem5.7} $\sR_0\setminus\sR^s_{0,\omega^0_c}=\widetilde{\mathcal{X}}_1\cup\cdots\cup\widetilde{\mathcal{X}}_r\cup\widetilde{\mathcal{Y}}_1\cup\cdots\cup\widetilde{\mathcal{Y}}_r$. In particular, $$\mathfrak{D}=\mathcal{X}_1+\cdots+\mathcal{X}_r+\mathcal{Y}_1+\cdots+\mathcal{Y}_r$$
where $\mathcal{X}_i:=\psi(\pi^{-1}(\widetilde{\mathcal{X}}_i)\cap\sR^s_{\omega^{0^+}} ),\, \mathcal{Y}_i:=\psi(\pi^{-1}(\widetilde{\mathcal{Y}}_i)\cap\sR^s_{\omega^{0^+}} )\,\subset\sU_{\mathbb{P}^1,\,\omega^{0^+}}$.
\end{lem}
\begin{proof} It is enough to show $\sR'_0\setminus\sR^s_{0,\omega^0_c}=\widetilde{\mathcal{X}}_1\cup\cdots\cup\widetilde{\mathcal{X}}_r\cup\widetilde{\mathcal{Y}}_1\cup\cdots\cup\widetilde{\mathcal{Y}}_{r-1}$, or equivalently
$\sR^s_{0,\omega^0_c}=\bigcap^r_{i=1}(\sR'_0\setminus\widetilde{\mathcal{X}}_i)\cap\bigcap^{r-1}_{i=1}(\sR'_0\setminus\widetilde{\mathcal{Y}}_i):=U$. Any point
$$q=(V\otimes\sO_{\mathbb{P}^1}(-N)\twoheadrightarrow E,\, Q_{\bullet}(E)_{y_1},\, Q_{\bullet}(E)_{z_1},\,Q_{\bullet}(E)_{z_2},\,Q_{\bullet}(E)_z)\in U$$
must satisfy the following condition
\begin{equation}\label{property linear algebra}\aligned &ev^{-1}_{y_1}(E_{y_1,1})\not\subset ev^{-1}_{z_1}(E_{z_1,r-i+1})\oplus ev^{-1}_{z_2}(E_{z_2,i})\subset{\rm H}^0(E)\\
&(1\le i\le r),\,\,ev^{-1}_{z_1}(E_{z_1,r-i})\cap ev^{-1}_{z_2}( E_{z_2,i})=0\,(1\leq i\leq r-1),\endaligned\end{equation}
which implies that the stabilizer ${\rm stab}(q)\subset {\rm PGL}(V)$ is trivial. Thus
$${\rm dim}\,{\rm PGL}(V)\cdot q={\rm dim}\,{\rm PGL}(V)={\rm dim}(U),\quad \forall\,\,q\in U.$$
It implies that $U={\rm PGL}(V)\cdot q$ ($\forall\,q\in U$). Take $q\in U\cap \sR^s_{0,\omega^0_c}$, we have
$$\sR^s_{0,\omega^0_c}={\rm dim}\,{\rm PGL}(V)\cdot q=U.$$
\end{proof}

To figure out the images $f(\mathcal{X}_i),\,\,f(\mathcal{Y}_i)\subset\mathrm{Flag}_{\vec n(z)}(E_z)$ of $f$, recall
$$f:\mathcal{U}_{\mathbb{P}^1,\,\omega^{0^+}}\rightarrow (\mathcal{U}_{
\mathbb{P}^1,\,\omega^0}\xrightarrow{\psi_q^{-1}})\mathrm{Flag}_{\vec n(z)}(E_z),$$
$q=(V\otimes\sO_{\mathbb{P}^1}(-N)\twoheadrightarrow E,\, Q_{\bullet}(E)_{y_1},\, Q_{\bullet}(E)_{z_1},\,Q_{\bullet}(E)_{z_2})\in \sR^s_{0,\omega^0_c}$ and let
\ga{5.13}{\bar{q}:=(E,\, Q_{\bullet}(E)_{y_1},\, Q_{\bullet}(E)_{z_1},\,Q_{\bullet}(E)_{z_2})\in \sU_{\mathbb{P}^1,\,\omega^0_c}}
be the unique isomorphic class of $\omega^0_c$-stable parabolic bundle, then the condition \eqref{property linear algebra} implies that
$L_i:=ev_{zz_1}^{-1}(E_{z_1,r-i})\cap ev_{zz_2}^{-1}(E_{z_2,i-1})\subset E_z$ ($1\le i\le r$) are one dimensional subspaces and
$$H_i:=ev^{-1}_{zy_1}(E_{y_1,1})\oplus \bigoplus_{1\le j\neq i,i+1\le r}L_j\subset E_z\quad (1\le i\le r-1)$$
are $r-1$ dimensional subspaces. On the other hand, for any injection $\varphi: \mathcal{O}_{\mathbb{P}^1}^{\oplus(r-2)}\oplus \mathcal{O}_{\mathbb{P}^1}(-1)\hookrightarrow E$
satisfying conditions
\ga{5.14}{E_{y_1,1}\subset {\rm Im}(\varphi_{y_1}),\quad E_{z_1,1}={\rm Im}(\varphi_{z_1}),\quad E_{z_2,1}={\rm Im}(\varphi_{z_2}),}
one can show that the $r-1$ dimensional subspace
\ga{5.15}{(E_{z_1,1}\cap E_{z_2,1}\subset) H_z:={\rm Im}(\varphi_{z})\subset E_z} is independent of the choices of $\varphi$. Then we have
\begin{lem}\label{lem5.8} Let $X_i:=f(\mathcal{X}_i),\,Y_i:=f(\mathcal{Y}_i)\subset\mathrm{Flag}_{\vec n(z)}(E_z)$. Then
\begin{itemize}
\item [(1)] $X_i=\{\, (\bar{q}, Q_{\bullet}(E)_z)\in \mathrm{Flag}_{\vec n(z)}(E_z)\,|\, E_{z,r-1}=L_i\,\}$ $(1\le i\le r)$.
\item[(2)] $Y_i=\{\, (\bar{q}, Q_{\bullet}(E)_z)\in \mathrm{Flag}_{\vec n(z)}(E_z)\,|\, E_{z,1}=H_i\,\}$ $(1\le i\le r-1)$,
$Y_r=\{\, (\bar{q}, Q_{\bullet}(E)_z)\in \mathrm{Flag}_{\vec n(z)}(E_z)\,|\, E_{z,1}=H_z\,\}$.
\end{itemize}
\end{lem}
\begin{proof} (1) Let $\mathcal{V}:=(\mathcal{V},Q_{\bullet}(\mathcal{V})_{y_1},\, Q_{\bullet}(\mathcal{V})_{z_1},\,Q_{\bullet}(\mathcal{V})_{z_2}, Q_{\bullet}(\mathcal{V})_{z})\in \mathcal{X}_i$ be a generic point. Then it is $\omega^{0^{+}}$-stable (but $(\mathcal{V},Q_{\bullet}(\mathcal{V})_{y_1},\, Q_{\bullet}(\mathcal{V})_{z_1},\,Q_{\bullet}(\mathcal{V})_{z_2})$ is not $\omega^0_c$-stable).
We will show, in the category of $\omega^0$-semistable parabolic bundles, $\mathcal{V}$ is $S$-equivalent to a parabolic bundle
$E:=(\bar{q}, Q_{\bullet}(E)_z)$ with $E_{z,r-1}=L_i$ (see \eqref{5.13} for the notion $\bar{q}$).

Let $\mathcal{V}'\cong\sO_{\mathbb{P}^1}^{\oplus r-1}\subset \sO_{\mathbb{P}^1}^{\oplus r}\cong\mathcal{V}$ be a sub-bundle of $\mathcal{V}$ such that
$${\rm H}^0(\mathcal{V}')=ev_{z_1}^{-1}(\mathcal{V}_{z_1,r-i+1})\oplus ev_{z_2}^{-1}(\mathcal{V}_{z_2,i})\subset {\rm H}^0(\mathcal{V}).$$
Then $ev^{-1}_{y_1}(\mathcal{V}_{y_1,1})\subset {\rm H}^0(\mathcal{V}')$ by the definition of $\mathcal{X}_i$. Thus
 ${n}_2^{\mathcal{V}'}(y_1)=1$, $$\aligned n^{\mathcal{V}'}_{r-i+2}(z_1)&=n^{\mathcal{V}'}_{r-i+3}(z_1)=\cdots =n^{\mathcal{V}'}_{r}(z_1)=1,\\
n^{\mathcal{V}'}_{i+1}(z_2)&=n^{\mathcal{V}'}_{i+2}(z_2)=\cdots =n^{\mathcal{V}'}_{r}(z_2)=1.
\endaligned$$
Since $2\leq \frac{\mathrm{par}\chi_{\omega^t}(\mathcal{V}')}{r-1}<\frac{\mathrm{par}\chi_{\omega^t}(\mathcal{V})}{r}=2+\frac{t}{r}$
holds for $t=0^+$, we must have $\frac{\mathrm{par}\chi_{\omega^t}(\mathcal{V}')}{r-1}=2$ and $n^{\mathcal{V}'}_{r-i+1}(z_1)=n_i^{\mathcal{V}'}(z_2)=n_r^{\mathcal{V}'}(z)=0$.
Then, in the category of $\omega^0$-semistable parabolic bundles, there is an exact sequence
$$0\rightarrow \mathcal{V}'\rightarrow \mathcal{V}\rightarrow (\mathcal{V}/\mathcal{V}'\cong\sO_{\mathbb{P}^1}, \{0,\frac{r-i}{r},\frac{i-1}{r}\}\cup \{\frac{1}{r}\}_{z})\rightarrow 0$$
with slopes ${\rm par}\mu_{\omega^0}(\mathcal{V}')={\rm par}\mu_{\omega^0}(\mathcal{V})={\rm par}\mu_{\omega^0}(\sO_{\mathbb{P}^1})$,
where $$\mathcal{V}':=(\mathcal{V}',Q_{\bullet}(\mathcal{V}')_{y_1},\, Q_{\bullet}(\mathcal{V}')_{z_1},\,Q_{\bullet}(\mathcal{V}')_{z_2}, Q_{\bullet}(\mathcal{V}')_{z})$$
is the parabolic sub-bundle of $\mathcal{V}$. The induced flag $0\subset \mathcal{V}'_{y_1,1}\subset \mathcal{V}'_{y_1}$,
$$0\subset \mathcal{V}'_{x,r'-1}\subset\cdots\subset \mathcal{V}'_{x,j}\subset \mathcal{V}'_{x,j-1}\subset\cdots\subset \mathcal{V}'_{x,1}\subset \mathcal{V}'_x,\quad x\in\{z_1,\,z_2,\,z\},$$ is of type $\vec n'(y_1)=(r'-1,1)$, $\vec n'(x)=(1,...,1)$ respectively ($r'=r-1$). We can check that flags in
$E':=(\mathcal{V}',Q_{\bullet}(\mathcal{V}')_{y_1},\, Q_{\bullet}(\mathcal{V}')_{z_1},\,Q_{\bullet}(\mathcal{V}')_{z_2})$
satisfy condition (\ref{property linear algebra}) (use $n^{\mathcal{V}'}_{r-i+1}(z_1)=n_i^{\mathcal{V}'}(z_2)=0$ and $\mathcal{V}\notin\mathcal{X}_j\cup\mathcal{Y}_j$ for $j\neq i$), which implies (see the proof of Lemma \ref{lem5.7}) that
$$(\mathcal{V}',Q_{\bullet}(\mathcal{V}')_{y_1},\, Q_{\bullet}(\mathcal{V}')_{z_1},\,Q_{\bullet}(\mathcal{V}')_{z_2})$$ is stable
respect to $\omega'_c:=(2r',\{\vec n'(x),\vec a'_c(x)\}_{x\in \{y_1,z_1,z_2\}})$. Then it is unique (up to isomorphism) by Proposition \ref{prop5.4}.
On the other hand, starting from $\bar{q}:=(E,\, Q_{\bullet}(E)_{y_1},\, Q_{\bullet}(E)_{z_1},\,Q_{\bullet}(E)_{z_2})\in \sU_{\mathbb{P}^1,\,\omega^0_c}$ and let
$\sO_{\mathbb{P}^1}\subset E$ be the sub-bundle such that
$$ev_{z_1}^{-1}(E_{z_1,r-i})\cap ev_{z_2}^{-1}(E_{z_2,i-1})={\rm H}^0(\sO_{\mathbb{P}^1})\subset {\rm H}^0(E)$$
(hence $n^{\sO_{\mathbb{P}^1}}_{r-i+1}(z_1)=n_i^{\sO_{\mathbb{P}^1}}(z_2)=1$, $n_2^{\sO_{\mathbb{P}^1}}(y_1)=0$), the flags in
$$\bar{q}':=(\bar{E}:=E/\sO_{\mathbb{P}^1}, Q_{\bullet}(\bar{E})_{y_1},\, Q_{\bullet}(\bar{E})_{z_1},\,Q_{\bullet}(\bar{E})_{z_2})$$
satisfy conditions (\ref{property linear algebra}). Hence $\bar{q}'$ is a $\omega'_c$-stable parabolic bundle, which must be isomorphic to
$(\mathcal{V}',Q_{\bullet}(\mathcal{V}')_{y_1},\, Q_{\bullet}(\mathcal{V}')_{z_1},\,Q_{\bullet}(\mathcal{V}')_{z_2})$. Then we can construct a flag
$Q_{\bullet}(E)_{z}$ of $E_z$ such that $(\bar{q}',Q_{\bullet}(\bar{E})_{z})\cong \mathcal{V}'$ where $Q_{\bullet}(\bar{E})_{z}$ is the image of $Q_{\bullet}(E)_{z}$
under $E_z\twoheadrightarrow \bar{E}_z:=(E/\sO_{\mathbb{P}^1})_z$. Thus
$$0\to (\sO_{\mathbb{P}^1}, \{0,\frac{r-i}{r},\frac{i-1}{r}\}\cup \{\frac{1}{r}\}_{z})\to (\bar{q},Q_{\bullet}(E)_{z})\to \mathcal{V}'\to 0$$
and $\mathcal{V}$ is $S$-equivalent to the parabolic bundle $(\bar{q},Q_{\bullet}(E)_{z})$ ($E_{z,r-1}=L_i$ since $n_r^{\mathcal{V}'}(z)=0$).
Conversely, given a generic $(\bar{q},Q_{\bullet}(E)_{z})$ with ${E_{z,r-1}=L_i}$, we can construct a $\mathcal{V}\in \mathcal{X}_i$ such
that $f(\mathcal{V})=(\bar{q},Q_{\bullet}(E)_{z})$. Then we are done, since $f$ is a projective morphism.

(2) The proof is similar with (1), we omit the detail here.
\end{proof}

Let $\widetilde{X}$ be the blowing-up of $X:=\mathrm{Flag}_{\vec n(z)}(E_z)$ along subvarieties $X_i,\,Y_i\,\subset X$ ($1\le i\le r$). Then there exist closed subsets
$Z_1\subset \widetilde{X}$ and $Z_2\subset \mathcal{U}_{\mathbb{P}^1,\,\omega^{0^+}}$ of codimensions at least two such that
$$\widetilde{X}\setminus Z_1\cong \mathcal{U}_{\mathbb{P}^1,\,\omega^{0^+}}\setminus Z_2.$$
Thus, to prove that $\sU_{\mathbb{P}^1,\,\omega^{0^+}}$ is split by a $(p-1)$-power, it is enough to show that $\widetilde{X}$ is split by a $(p-1)$-power. To simplify notation, let $W:={\rm H}^0(E)$ and identify $W$ with $E_z$ (through $ev_z:W\cong E_z$). Then
$$X=\mathrm{Flag}_{(1,...,1)}(W)=\{\,\upsilon=(0\subset V_1\subset\cdots\subset V_{r-1}\subset V_r=W)\,\}.$$
Let $L_i\subset W$ ($1\leq i\leq r$) be subspaces of dimension one such that ${W=L_1+\cdots+ L_r}$,
$L\subset W$ be a generic subspace of dimension one and $H\subset W$ be a subspace of dimension $r-1$ such that $L,\,L_1,\,L_r\nsubseteq H$ and $L_i\subset H$ for $i\neq 1,\,r$. Let $Y_r:=\{\,\upsilon\in X\,|\, V_{r-1}=H\,\}\subset X,$
\begin{equation}\label{5.16}\aligned &X_i:=\{\,\upsilon\in X\,|\, V_1=L_i\,\}\subset X\,\,\,\,\,\,\,\, (1\le i\le r),\\
&Y_i:=\{\,\upsilon\in X\,|\, V_{r-1}=H_i\,\}\subset X \,\,(1\le i\le r-1)\endaligned\end{equation}
where $H_i=L+L_1+\cdots+L_{i-1}+L_{i+2}+\cdots+L_r$.

In order to show the blowing-up $\widetilde{X}$ of $X$ along $X_i,\,Y_i\subset X$ ($1\le i\le r$) is split by a $(p-1)$-power, we have to
construct a section $\sigma\in {\rm H}^0(X,\omega_X^{-1})$ such that (1) its vanishing orders at $X_i$, $Y_i$ ($1\leq i\leq r$) satisfy
\begin{equation}\label{5.17}\aligned
\mathrm{ord}_{X_i}(\sigma)&\geq r-2=(r-1)-1\\
\mathrm{ord}_{Y_i}(\sigma)&\geq r-2=(r-1)-1
\endaligned\end{equation}
where $r-1$ is the codimension of $X_i$, $Y_i\subset X$; and (2) $\sigma^{p-1}$ splits $X$ (see \cite[1.3.E Exercises (4*)]{M.S}).


To construct $\sigma$, let $0\subset \mathcal{V}_1\subset\cdots\subset \mathcal{V}_{r-1}\subset \mathcal{V}_r=W\otimes \sO_X$ be the universal
flag on $X=\mathrm{Flag}_{(1,\ldots,1)}(W)$, and for $1\leq i\leq \frac{r}{2}$, we set
\begin{equation}\label{5.18}\aligned W_{r+1-2i}:=\sum_{ j\in  [i+1,r+1-i]} L_j,\, &W_{r-2i}:=\sum_{ j\in [i+1,r-i] }L_j,\,\,
\\ S_{2i-1}:=L+\sum_{j\notin [i,r+1-i]} L_j,\, &S_{2i}:=L+\sum_{j\notin [i+1,r+1-i]}L_j.\endaligned\end{equation}
Define $D_i\subset X$ to be the locus where
$$W_i\otimes\sO_X\hookrightarrow W\otimes\sO_X\rightarrow W\otimes\sO_X/\sV_{r-i}$$ is not an isomorphism, $E_i\subset X$ to be the locus where
$$S_i\otimes\sO_X\hookrightarrow W\otimes\sO_X\rightarrow W\otimes\sO_X/\sV_{r-i}$$ is not an isomorphism. Then $\sO_X(D_i)\cong{\rm det}(\sV_{r-i})^{-1}\cong\sO_X(E_i)$ and
$$\omega^{-1}_X=\bigotimes_{i=1}^{r-1}{\rm det}(\sV_i)^{-2}=\bigotimes^{r-1}_{i=1}\sO_X(D_i)\otimes \bigotimes^{r-1}_{i=1}\sO_X(E_i). $$
Let $d_i$ (resp. $e_i$) be canonical section of $D_i$ (resp. $E_i$), then $$\sigma:=d_1\cdots d_{r-1}\cdot e_1\cdots e_{r-1}\in {\rm H}^0(\omega^{-1}_X)$$
satisfies the condition \eqref{5.17}. Indeed, for $1\leq i\leq \frac{r}{2}$,
$$\aligned &[i,r+1-i]-[i+1,r+1-i]=\{i\},\,\,  \\& [i+1,r+1-i]-[i+1,r-i]=\{r+1-i\},\endaligned$$
thus for any fixed $1\leq j\leq r$, there is at most one $k_0\in [1,r-1]$ such that $L_j\nsubseteq W_{k_0}\cup S_{r-k_0}$ (see (\ref{5.18})).
Choose $k\neq k_0\in  [1,r-1]$, then $d_k\cdot e_{r-k}$ vanishes along $X_j$ (see (\ref{5.16})) and it is clear that $\mathrm{ord}_{X_j}(\sigma)\geq r-2$. One can also check that there are at least $r-2$ indexes $k\in [1,r-1]$ such that
$W_k\subset H_j$ or $S_{r-k}\subset H_j$ for any fixed $1\leq j\leq r-1$. 
All together, it is enough to show
\begin{lem}\label{lem5.9}
The section $\sigma^{p-1}\in {\rm H}^0(\omega^{1-p}_X)$
splits $X=\mathrm{Flag}_{(1,\ldots, 1)}(W)$.
\end{lem}
\begin{proof}
To prove the lemma, we will use repeatedly the fact: let $X$ be a normal projective variety and $Y\subset X$ be a normal divisor with canonical section $s\in {\rm H}^0(\sO_X(Y))$. If $\tilde{\alpha}\in {\rm H}^0(\omega^{-1}_X(-Y))$ is a lift of $\alpha\in{\rm H}^0(\omega^{-1}_Y)$ and $\alpha^{p-1}$ splits $Y$, then $(\tilde{\alpha}\cdot s)^{p-1}$ splits $X$.

Let $w=(0\subset W_1\subset \cdots \subset W_{r-1}\subset W_r=W)$ (see (\ref{5.18})) be the fixed point of $\mathrm{Flag}_{(1,\ldots, 1)}(W)$ and, for $1\le k\le r$, let
$$R_k=\{\,v\in \mathrm{Flag}_{(1,\ldots,1)}(W)\,|\, V_i=W_i,\,\,k\le i\le r\,\},$$
then $\mathrm{Flag}_{(1,\ldots, 1)}(W)=R_r\supset R_{r-1}\supset\cdots\supset R_2\supset R_1=\{w\}$. Consider
\ga{5.19}{R_k=Z_{k,0}\supset Z_{k,1}\supset Z_{k,2}\supset\cdots\supset Z_{k,k-2}\supset Z_{k,k-1}=R_{k-1}}
where $Z_{k,i}=\{\,v\in R_k\,|\, V_i\subset W_{k-1}\,\}$ ($0\le i<k$). On $R_k$, we have
\ga{5.20}{\sV_{i+1}\hookrightarrow \sV_k=W_k\otimes \sO_X\twoheadrightarrow W_k/W_{k-1}\otimes \sO_X \,\,(0\le i<k)}
 which, on $Z_{k,i}$, induces ${(\sV_{i+1}/\sV_i)|_{Z_{k,i}}\to (W_k/W_{k-1}\otimes \sO_X)|_{Z_{k,i}}}$ whose zero locus is
$Z_{k,i+1}\subset Z_{k,i}$. Let ${\sL_{k,i}=(\sV_{i+1}/\sV_i)^{-1}\otimes(\sW_k/W_{k-1}\otimes \sO_X)}$,
$$\sL_{k,i}|_{Z_{k,i}}=\sO_{Z_{k,i}}(Z_{k,i+1}), \quad 1_{Z_{k,i+1}}\in {\rm H}^0(\sO_{Z_{k,i}}(Z_{k,i+1}))$$
where $1_{Z_{k,i+1}}\in {\rm H}^0(\sO_{Z_{k,i}}(Z_{k,i+1}))$ is the canonical section of $Z_{k,i+1}$. Let
$$\sM_{k,i}={\rm det}(\sV_i)^{-1}\otimes{\rm det}(W_{k}/W_{k-i}\otimes \sO_X)\quad (1\le i\le k-1),$$
which, on $R_k$, have global section $1_{C_{k,i}}\in {\rm H}^0(\sM_{k,i}|_{R_k})={\rm H}^0(\sO_{R_k}(C_{k,i}))$ where
$C_{k,i}$ is the divisor of $1_{C_{k,i}}$. By easy computations, we have
$$\aligned &\omega^{-1}_{R_k}(-C_{k,1}-\cdots-C_{k,k-1})|_{R_{k-1}}=(\omega^{-1}_{R_k}\otimes\sL^{-1}_{k,0}\otimes\cdots\otimes\sL^{-1}_{k,k-2})|_{R_{k-1}}
\otimes\\ &\sO_{R_{k-1}}(-C_{k-1,1}-\cdots-C_{k-1,k-2})\quad (1<k\le r)\endaligned $$
and $\omega^{-1}_{R_k}=(\omega^{-1}_{R_{k+1}}\otimes \sL^{-1}_{k+1,0}\otimes\cdots\otimes\sL^{-1}_{k+1,k-1})|_{R_k}$. Thus
\ga{5.21}{\omega^{-1}_{R_r}(-C_{r,1}-\cdots-C_{r,r-1})|_{R_k}=\omega^{-1}_{R_k}(-C_{k,1}-\cdots-C_{k,k-1}).}
Let $\delta_k=\tilde{\sigma}|_{R_k}\cdot 1_{C_{k,1}}\cdot1_{C_{k,2}} \cdots 1_{C_{k,k-1}} \in {\rm H}^0(\omega^{-1}_{R_k})$ ($1\le k\le r$) where $\widetilde{\sigma}=e_1\cdots e_{r-1}$. Note $\delta_1=\widetilde{\sigma}|_{\{w\}}$ which is not zero by definition of $e_i$ ($1\leq i\leq r-1$) and $\delta_r=\sigma$, it is enough to show that $\delta^{p-1}_k$ splits $R_k$ if $\delta^{p-1}_{k-1}$ splits $R_{k-1}$ for any $k>1$. To prove that $\delta^{p-1}_k$ splits $R_k$, use filtration \eqref{5.19}, note $C_{k,1}=Z_{k,1}$ and $\delta_{k,0}:=\delta_k=\tilde{\delta}_{k,1}\cdot 1_{Z_{k,1}}$ where $$\tilde{\delta}_{k,1}:=\tilde{\sigma}|_{Z_{k,0}}\cdot 1_{C_{k,2}} \cdots 1_{C_{k,k-1}}\in {\rm H}^0(\omega^{-1}_{Z_{k,0}}(-Z_{k,1})),$$ it is enough to prove that $\delta^{p-1}_{k,1}$ splits $Z_{k,1}$ where $$\delta_{k,1}:=\tilde{\delta}_{k,1}|_{Z_{k,1}}\in{\rm H}^0(\omega^{-1}_{Z_{k,1}}).$$
We will show that there exists $\tilde{\delta}_{k,2}\in {\rm H}^0(\omega^{-1}_{Z_{k,1}}(-Z_{k,2}))$ such that
$$\delta_{k,1}=\tilde{\delta}_{k,2}\cdot 1_{Z_{k,2}}.$$
Indeed, on $Z_{k,i}$ ($1\le i<k-1$), we have the following diagram
$${\small \begin{CD}
0&@>>>& \sV_i &@>>> &\sV_{i+1} &@>>>&\sV_{i+1}/\sV_i &@>>>&0\\
&&&&@VVV&&@V\phi_{k,i+1}VV&&@VVV\\
0&@>>>&\frac{W_{k-1}}{W_{k-1-i}}\otimes \sO_X &@>>>&\frac{W_k}{W_{k-1-i}}\otimes\sO_X&@>>>&{\Tiny \frac{W_k}{W_{k-1}}\otimes \sO_X}&@>>>&0
\end{CD}}$$
where $\sV_i\to W_{k-1}/W_{k-1-i}\otimes \sO_X$ (on $Z_{k,i}$) is induced by
 $$\phi_{k,i+1}:\sV_{i+1}\hookrightarrow \sV_k=W_k\otimes\sO_X\twoheadrightarrow W_k/W_{k-1-i}\otimes \sO_X,$$
thus there is a global section $s_{k-1,i}\in {\rm H}^0(\sM_{k-1,i}|_{Z_{k,i}})$ such that
\ga{5.22}{1_{C_{k,i+1}}|_{Z_{k,i}}=s_{k-1,i}\cdot 1_{Z_{k,i+1}},\quad s_{k-1,i}|_{R_{k-1}}=1_{C_{k-1,i}}.}
In particular, $\delta_{k,1}:=\tilde{\delta}_{k,1}|_{Z_{k,1}}=(\tilde{\sigma}|_{Z_{k,0}}\cdot 1_{C_{k,3}} \cdots 1_{C_{k,k-1}})|_{Z_{k,1}}\cdot 1_{C_{k,2}}|_{Z_{k,1}}=
\tilde{\delta}_{k,2}\cdot 1_{Z_{k,2}}$ where
$\tilde{\delta}_{k,2}:=(\tilde{\sigma}|_{Z_{k,0}}\cdot 1_{C_{k,3}} \cdots 1_{C_{k,k-1}})|_{Z_{k,1}}\cdot s_{k-1,1}.$
In fact, let
$$\tilde{\delta}_{k,i+1}:=(\tilde{\sigma}|_{Z_{k,0}}\cdot 1_{C_{k,i+2}} \cdots 1_{C_{k,k-1}})|_{Z_{k,i}}\cdot (s_{k-1,1}\cdots s_{k-1,i})|_{Z_{k,i}}$$
and $\delta_{k,i+1}:=\tilde{\delta}_{k,i+1}|_{Z_{k,i+1}}$ ($0\le i\le k-1$), by the fact stated at beginning of the proof, we know that $\delta_{k,i}=\tilde{\delta}_{k,i+1}\cdot 1_{Z_{k,i+1}}$ provides a $(p-1)$-power splitting of $Z_{k,i}$ if $\delta_{k,i+1}=\tilde{\delta}_{k,i+1}|_{Z_{k,i+1}}$ provides a $(p-1)$-power splitting of $Z_{k,i+1}$.
Since $s_{k-1,k-1}|_{R_{k-1}}$ is the constant section of $\sM_{k-1,k-1}|_{Z_{k,k-1}}={\rm det}(\sV_{k-1})^{-1}\otimes{\rm det}(W_{k-1}\otimes\sO_X)|_{R_{k-1}}$, we are done by
$$\aligned\delta_{k,k-1}&=\tilde{\sigma}|_{Z_{k,k-1}}\cdot(s_{k-1,1}\cdots s_{k-1,k-2}\cdot s_{k-1,k-1})|_{Z_{k,k-1}}\\
&=\tilde{\sigma}|_{R_{k-1}}\cdot 1_{C_{k-1,1}}\cdot\cdots \cdot 1_{C_{k-1,k-2}}=\delta_{k-1}.\endaligned $$
\end{proof}

\begin{proof}[Proof of Theorem \ref{thm5.1}]
By Hecke transformation, we can assume $d=0$. Then it is enough to show that $\sU_{\mathbb{P}^1,\,\omega^{0^+}}$ is split by a $(p-1)$-power (see Lemma \ref{lem5.6}),
which is equivalent to prove, by Lemma \ref{lem5.8}, the blowing-up of $X=\mathrm{Flag}_{(1,\ldots, 1)}(W)$ along $X_i,\,Y_i\subset X$ ($1\le i\le r$) is split by a $(p-1)$-power. Then we are done by Lemma \ref{lem5.9}.
\end{proof}

\bibliographystyle{plain}

\renewcommand\refname{References}

\end{document}